\documentclass{article}[11pt]
\usepackage{amsmath, latexsym, amsfonts, amssymb, amsthm, amscd, color, bbm}
\usepackage{psfrag,graphicx, appendix}

\newtheorem {lemme} {Lemma} [section]
\newtheorem {theoreme} {Theorem} [section]
\newtheorem {proposition} {Proposition} [section]

\newtheorem{remark}{Remark}[section] 
\newtheorem {definition} {Definition} [section]

\newcommand{\E}{\mathbb {E}}

\newcommand{\R}{\mathbb {R}}
\newcommand{\C}{\mathbb {C}}

\newcommand{\vers}{\mathop{\longrightarrow }}
\newcommand{\bl}{\begin{lemme}}
\newcommand{\el}{\end{lemme}}

\newcommand{\brem}{\begin{remark}}
\newcommand{\erem}{\end{remark}}
\newcommand{\bconj}{\begin{conjecture}}
\newcommand{\econj}{\end{conjecture}}
\newcommand{\bdefi}{\begin{definition}}
\newcommand{\edefi}{\end{definition}}
\newcommand{\bt}{\begin{theoreme}}
\newcommand{\bfa}{\begin{fact}}
\newcommand{\efa}{\end{fact}}

\newcommand{\et}{\end{theoreme}}
\newcommand{\bp}{\begin{proposition}}
\newcommand{\ep}{\end{proposition}}
\newcommand{\be}{\begin{equation}}

\newcommand{\ee}{\end{equation}}

\renewcommand{\Re}{\mathrm{Re}}
\renewcommand{\Im}{\mathrm{Im}}

\title{Fluctuations at the edges of the spectrum of the full rank deformed GUE}
\author{Mireille Capitaine\footnote{CNRS, Institut de Math\'ematiques de Toulouse, Universit\'e Paul Sabatier
118, Route de Narbonne
F-31062 Toulouse Cedex 9, mireille.capitaine@math.univ-toulouse.fr 
}  and  Sandrine P\'ech\'e\footnote{Universit\'e Paris Diderot, LPMA, 5 Rue Thomas Mann, 75013 Paris, France, sandrine.peche@math.univ-paris-diderot.fr}}
\begin{document}
\maketitle
\begin{abstract}
We consider a full rank deformation of the GUE $W_N+A_N$ where $A_N$ is a full rank Hermitian matrix of size $N$ and $W_N$ is a GUE. The empirical eigenvalue distribution $\mu_{A_N}$ of $A_N$ converges to a probability distribution $\nu$. We identify all the possible limiting eigenvalue statistics at the edges of the spectrum, including  outliers, edges and merging points of connected components of the limiting spectrum. The results are stated in terms of a deterministic equivalent of the empirical eigenvalue distribution of $W_N+A_N$, namely the free convolution of the semi-circle distribution and the empirical eigenvalues distribution of $A_N$.
\end{abstract}
\section{Introduction and results}
\subsection{Motivations}
Enormous progress has been accomplished in the very recent years in the study of asymptotic spectral properties of large random matrices. A Hermitian Wigner random matrix is a $N\times N$ matrix $W_N=\frac{1}{\sqrt N}(W_{ij})_{i,j=1}^N$,with i.i.d. entries off the diagonal $W_{ij}, i<j$ (modulo the symmetry assumption) and independent diagonal real entries. 
The entries are standardized to be centered and of variance $\sigma^2$.
The asymptotic local properties of the spectrum of Wigner random matrices are now quite well understood thanks to the fantastic work of Erd\"os-Schlein-Yau (see \cite{Yauetal1},\cite{Yauetal} and references therein) and Tao-Vu \cite{TaoVu}.  In particular, it is known (assuming that the matrix elements admit enough moments) that the fluctuations of eigenvalues in the bulk or at the edges of the spectrum are  universal. In particular, they coincide with those identified for a Gaussian (GUE) matrix with variance $\sigma^2$. In other words, the limiting asymptotic spectral properties of a Wigner matrix in the large $N$ limit do not depend on the detail of the distribution of the matrix elements $W_{ij}$, $1\le i,j \le N.$ \\
In this article, we are interested in deformed random matrix ensembles.  
A \emph{deformation } of a standard random matrix can be more or less understood as the modification of the distribution of some of the entries of a Wigner matrix. The set of possible deformations is non exhaustive (one can force some of the entries to be zero such as for sparse matrices) but we here restrict to some additive deformations. 
More precisely, we consider a matrix $A_N$ of size $N$, which is deterministic. Our study could be extended to the case where it is random but we do not wish to pursue this direction here.  We consider the deformed matrices $$W_N+A_N ,$$ where $W_N$ is a standard Wigner matrix.  The question is to understand the asymptotic properties of the eigenvalues and eigenvectors of the deformed matrix, knowing that of $A_N$ and $W_N$. Such ensembles have first been introduced by \cite{BrezinHikami}, and \cite{Johansson} when $W_N$ is a GUE.

\paragraph{}In the case where $A_N$ is a fixed rank (independent of the size $N$) matrix, the asymptotic properties of the spectrum are quite clear. Finite rank perturbed ensembles have first been considered in \cite{BBP} (see also \cite{BaikSilverstein} and \cite{Peche}). First, the global properties of the spectrum are not impacted by $A_N$. Indeed, denoting by $\lambda_1\geq \lambda_2 \geq \cdots \geq \lambda_N$ the ordered eigenvalues of $W_N +A_N$, the empirical eigenvalue distribution $\mu_N:=\frac{1}{N}\sum_{i=1}^N \delta_{\lambda_i}$ still converges (as in the case where $A_N=0$) to the semi-circle distribution with density $\sigma_{sc}(x)=\frac{1}{2\pi \sigma^2}\sqrt{4\sigma^2-x^2} \mathbbm{1}_{|x|\leq 2\sigma}.$ 
The asymptotic local eigenvalue statistics of eigenvalues in the bulk of the spectrum are also unchanged by the deformation matrix $A_N$.
Only the local behavior of the spectrum at the edges may be impacted by the deformation $A_N$, as we now explain.
The deformation $A_N$ may cause some eigenvalues to separate from the bulk of the spectrum. Each eigenvalue of $A_N$ greater than $\sigma$ is called a \emph{spike}. To each spike $\theta_i$ of $A$ such that $|\theta_i|>\sigma$ (if it exists) there corresponds an eigenvalue $\lambda_i $ satisfying
$$\lambda_i \to (\sigma+\dfrac{\sigma^2}{\theta_i})$$ a.s. Such eigenvalues $\lambda_i$ outside the support of the semi-circle distribution are called \emph{outliers}.
Interestingly, \cite{CapDonFeral} and then  \cite{RenfrewSos1}, \cite{RenfrewSos2} have proved that the fluctuations of spikes are not universal in general.  More precisely 
$$\sqrt N \left( \lambda_i -(\sigma+\dfrac{\sigma^2}{\theta_i})\right )\stackrel{d}{\to}\mu,$$
where the distribution $\mu$ may depend explicitly on the distribution of the matrix elements $W_{ij}$. It can be shown that eigenvectors of the matrix $A_N$ play a fundamental role in the universality/non universality of the deformation matrix $A_N$.
On the contrary, when there is no spike, the limiting distribution of extreme eigenvalues is the same as in the non deformed case. In particular, extreme eigenvalues stick to the bulk of the spectrum. The scale of their fluctuations is $N^{-2/3}$ and the limiting distribution of the largest (and smallest) eigenvalues is the Tracy-Widom distribution, provided the matrix elements $W_{ij}$ admit enough moments.
A complete study of such deformed ensembles has been achieved in \cite{KnowlesYin} and \cite{KnowlesYin2} and we refer the reader to these articles for a complete state of the art in finite rank deformations of Wigner matrices.

\paragraph{}The study of deformed ensembles extends to the case where the matrix $A_N$ has low rank $r_N <<N, r_N \to \infty$ (see \cite{Peche} e.g.) or full rank i.e. when $r_N=O(N)$. In this case, it is natural to assume that the empirical eigenvalue distribution of $A_N$ has a weak limit as $N \to \infty$, which is possibly $\delta_0$. Denote by $y_1\geq y_2\geq \cdots \geq y_N$ the ordered eigenvalues of $A$. 
Let $\mu_{A_N}=\frac{1}{N}\sum_{i=1}^N \delta_{y_i}.$ We assume the norms of $ (A_N)_{N}$ are uniformly bounded and  that there exists a probability distribution $\nu$ on $\mathbb{R}$ such that 
$$\mu_{A_N}\underset{N \to \infty}{\overset{w}{\to}}\nu.$$ Let us diagonalize $A_N$ through $A_N=V \text{diag}(y_1, \ldots, y_N) V^*$.
Roughly speaking the deformed model is now understood in the sense that $A_N$ is a "small" perturbation of the matrix 
$W_N+VA_0V^*$ where $A_0$ would be a diagonal matrix made up with quantiles of the probability $\nu.$
The asymptotic global behavior of the spectrum is well-known in this case. Indeed, let $\mu_{N}$ be the empirical eigenvalue distribution of $W_N+A_N.$ Its Stieltjes transform is
$$m_N(z):= \int \frac{1}{z-y}d \mu_{N}(y), \Im z \not=0.$$ According to  \cite{Voiculescu91, AGZ09}, $m_N$ converges as $N\to \infty$ to the Stieltjes transform $m_{\tau}$ of a probability distribution $\tau$, called the free convolution of $\nu$ and the semi-circle distribution. This probability distribution $\tau$ is uniquely characterized by a fixed point equation satisfied by $m_{\tau}$, as we review in Section \ref{sec: freeconvgauss}; it has a density $p$. We emphasize that the support of the probability distribution $\tau$ may have distinct connected components, depending on $\nu.$ 

\paragraph{}The question of the asymptotic behavior 
of extreme eigenvalues naturally arises in this setting also. This question has been much less investigated actually. So far, only the case where $W_N$ is a GUE has been investigated. \\
In \cite{Shcherbina2}, the author considers the case where $\mu_{A_N}$ concentrate quite fast to the measure $\nu$. In particular, there are no spikes. When $W_N$ is a GUE, she investigates the local edge regime which deals with the behavior of the eigenvalues near any extremity point $u_0$ of a connected component of $\text{supp}(\tau)$. More precisely let some $\epsilon >0$ be given and assume that   either  \begin{eqnarray}& & p(u)>0, ~~\forall u \in ]u_0; u_0+ \epsilon[, ~~\text{and}~~p(u)=0, ~~ \forall  u \in ]u_0-\epsilon; u_0],\label{re}\\
&\text{or }&p(u)>0, ~~\forall u \in ]u_0-\epsilon; u_0[, ~~\text{and}~~p(u)=0, ~~ \forall u \in [u_0; u_0+ \epsilon[.\end{eqnarray}
\cite{Shcherbina2} makes a technical assumption on the uniform convergence of the Stieltjes transform of $\mu_{A_N}$ to $m_{\nu}$: 
\begin{equation}\label{Scherbina}\sup_{z \in K} |m_{\mu_{A_N}}(z)-m_{\nu }(z)|\leq N^{-2/3-\epsilon}, \end{equation} where $K$ is some compact subset of the complex plane at a positive distance of the support of $\nu.$ This is a rather strong assumption on the rate of convergence of $\mu_{A_N}$ to $\nu$.
\cite{Shcherbina2} proves  that the joint distribution of the largest (or smallest) eigenvalues converging to $u_0$ have universal asymptotic behavior, characterized by the famous Tracy-Widom distribution. We note that \cite{Shcherbina} also investigates the asymptotic spacing distribution of eigenvalues in the bulk of the spectrum. The same behavior as for non deformed ensemble is obtained (and described by the sine kernel). The extension to a non Gaussian matrix $W$ has recently been obtained by \cite{RenfrewVu} in the case where $A_N$ is diagonal.\\
In \cite{Adleretal1} and \cite{Adleretal2}, the authors consider the case where $\mu_{A_N}=\nu$ is a finite combination of Dirac delta masses. They identify different possible limiting statistics at the edges of the support of $\tau$, after suitable normalization of the eigenvalues. If $u_0$ is a point such that $p(u)=0, \: u_0-\epsilon \le u\le u_0 $, $p(u)>0, u_0<u\le u_0+\epsilon$ for some $\epsilon >0$, the asymptotic distribution of eigenvalues close to $u_0$ is the Tracy-Widom distribution. The authors also consider the case where $u_0$ is a point where two connected components of $\text{supp}(\tau)$ merge so that 
$p(u)>0 , \forall u \in (u_0-\epsilon, u_0+\epsilon)\setminus \{u_0\} $ and $p(u_0)=0$. In this case, the limiting eigenvalue statistics are described by the so-called Pearcey kernel (whose definition is reviewed hereafter).

\paragraph{}
In both cases, a strong assumption is made on the rate of convergence of $\mu_{A_N}$ to $\nu$. We here remove this assumption. We identify all the possible limiting eigenvalue statistics at the edges of the spectrum of the deformed GUE, namely at a spike, at the edge of a connected component of the support or at a point where two connected components merge. We emphasize that we do not make any assumptions on the rate of convergence of  $\mu_{A_N}$ to $\nu$.
To state our results, we use a deterministic equivalent of the empirical eigenvalue distribution of $M_N$. This equivalent is the free convolution of the semi-circle distribution and $\mu_{A_N}$.\\
The choice of the deformed GUE is motivated by the fact that all eigenvalues statistics can be explicitly computed for this ensemble of deformed random matrices. We expect that one can extend these results to full rank deformations of an arbitrary Wigner matrix, as in the fixed rank case (with universal or non universal results). We intend to consider this general case in a forthcoming paper. The techniques needed are completely different.

\subsection{Model and results}\label{model} 
We consider the 
following deformed GUE ensemble $$M_N=X_N+ A_N,$$ where 
\begin{itemize}
\item[($H_1$)] $X_N=\frac{1}{\sqrt{N}} W_N $ where $W_N$ is a $N \times N$ GUE matrix: the random variables
$(W_N)_{ii}$, $\sqrt{2} (\Re (W_N)_{ij})_{i < j}$, $\sqrt{2}
(\Im (W_N)_{ij})_{i<j}$
are  i.i.d., with gaussian distribution of variance $1$ and mean 0.
\item[($H_2$)] $A_N$ is a deterministic Hermitian matrix whose eigenvalues $y_i=y_i(N)$, $1\leq i\leq N$,
are such that the spectral measure $\mu _{A_N} := \frac{1}{N} \sum_{i=1}^N \delta _{y_i}$ 
converges weakly to some probability measure $\nu $ 
with compact support. We assume that
\begin{equation}\label{suppnu} \forall t \in  \text{ supp }(\nu), ~~\lim_{\epsilon \to 0}\int \frac{d\nu(x)}{(t+i\epsilon-x)^2}>1,\end{equation}
where ${\rm supp}(\nu )$ denotes the support of $\nu $. 
\item[($H_3$)]  We also assume that there exists a fixed integer $r\geq 0$ (independent from $N$) and an integer $0\leq J\leq r$ such that the following holds.  There are $J$ fixed real numbers $\theta _1 > \ldots > \theta _J$ 
independent of $N$ which are outside the support of $\nu $ and such that each $\theta _j$ 
is an eigenvalue of $A_N$ with a fixed multiplicity $k_j$ (with $\sum_{j=1}^J k_j=r$). 
The $\theta_j$'s are called the spikes or the spiked eigenvalues of $A_N$ and we set $$\Theta =\{ \theta _j, \, 1 \leq j \leq J \}. $$
The remaining $N-r$ eigenvalues of $A_N$, denoted by $\beta _j(N)$, $j=1, \ldots, N-r$, satisfy 
$$\max _{1\leq j\leq N-r} {\rm dist}(\beta _j(N),{\rm supp}(\nu ))\vers _{N \rightarrow \infty } 0.$$
\end{itemize} 
Denote by  $\mu _{sc}$ the semicircle distribution whose density is given by
\begin{equation}\label{scl}
\frac{d\mu _{sc }}{dx}(x)= \frac{1}{2 \pi } \sqrt{4- x^2} 
\, 1 \hspace{-.20cm}1_{[-2 , 2]}(x).
\end{equation}
According to \cite{AGZ09}, the spectral distribution of $M_N$ 
weakly converges  almost surely to the so-called \emph{free convolution} $\mu _{sc} \boxplus \nu $ which has a continuous density  $p$ (see \cite{Biane97b}).
We recall some important facts about the free convolution with a semi-circular distribution  in Section \ref{sec: freeconvgauss}.

 \paragraph{}

We are now in position to state our results.
Let first consider a real number $d$ which is a right edge of $\text{supp}(\mu_{sc}\boxplus \nu)$ that is which satisfies (\ref{re}).
Assume moreover that for any $\theta_j$ such that $\int \frac{d\nu(s)}{(\theta_i -s)^2} =1$, we have $d \neq \theta_j + m_\nu (\theta_j)$.
 We show in Proposition \ref{rightedge} that 
for $\eta$ small enough, for all large $N$, there exists a  unique right edge $d_N$ of  $\text{supp}(\mu_{sc}\boxplus \mu_{A_N})$ in
$]d-\eta; d+\eta[$.
We derive the asymptotic distribution of eigenvalues in the vicinity of $d_N$.
Before exposing our results, we need a few notations. Let $Ai(u)$ be the Airy function defined by
\begin{equation}\label{eq:Airyintrep}
  Ai(u) = \frac{1}{2\pi } \int e^{iua+i\frac{1}3 a^{3}}da
\end{equation}
where the contour is from $\infty e^{5i\pi/6}$ to $\infty
e^{i\pi/6}$. The \emph{Airy kernel} (see e.g.
\cite{TW}) is then given by
\begin{equation}
  \mathbf{A}(u,v) = \frac{Ai(u)Ai'(v)-Ai'(u)Ai(v)}{u-v}=\int_0^\infty Ai(u+z)Ai(z+v) dz.
\end{equation}
Let $\mathbf{A}_x$ be the operator acting on $L^2((x,\infty))$
with kernel $\mathbf{A}(u,v)$. 
The GUE Tracy-Widom distribution for the largest eigenvalue is 
(\cite{TW})
\begin{equation}\label{eq:TWGUE}
  F_0(x)= \det(1-\mathbf{A}_x)=F_{GUE}(x).
\end{equation}
We refer to \cite{TW} for the more complicated definition of the GUE distribution for the $k$ largest eigenvalues ($k>1$).

We first prove the following result. Let $k$ be a given fixed integer. Let $\lambda_{max}\geq \lambda_{max-1}\geq \cdots \lambda_{max-k+1}$ denote the $k$ largest of those eigenvalues of $M_N$  converging to $d.$
\bt \label{theo: A}There exists $\alpha >0$ depending on $d_N$ only such that 
 the vector $$\frac{N^{2/3}}{\alpha} \left( \lambda_{max}-d_N, \lambda_{max-1}-d_N , \ldots, \lambda_{max-k+1}-d_N  \right)$$ converges in distribution as $N \to \infty$ to the so-called Tracy-Widom GUE distribution for the $k$ largest eigenvalues.
\et 
\brem Condition \ref{suppnu} is necessary to obtain Tracy-Widom asymptotics at the edges of the spectrum. If condition \ref{suppnu} fails e.g. at the top edge of the spectrum, meaning that the density of $\nu$ vanishes too fast at the edge, the limiting eigenvalue statistics at the edge can be proved to be Gaussian. 
\erem

We now turn to the behavior of outliers. Let $\theta_{i}$ be a spiked eigenvalue with multiplicity $k_i$, such that 
$\int \frac{1}{(\theta_i-x)^2}d\nu(x) <1$. In \cite{CDFF}, the authors prove that the spectrum of $M_N$ exhibits $k_i$ eigenvalues in a neighborhood of \be \label{defrho} \rho_{\theta_{i}}=\theta_i+\int\frac{d\nu(x)}{\theta_i-x}.\ee
Note that such a result is obtained when the support of $\nu$ has a finite number of connected components. However this assumption can be easily relaxed (see Remark \ref{extension}).
In Proposition \ref{outlier}, we prove that for $\epsilon>0$ small enough, for all large $N$, $\text{supp}(\mu_{sc}\boxplus\mu_{A_N})$ has a unique connected component $[L_i(N); D_i(N)]$ inside $]\rho_{\theta_i} -\epsilon;  \rho_{\theta_i} +\epsilon[$.
Define 
\be \label{defrhoN} \rho_N(\theta_{i})=z+\frac{1}{N}\sum_{y_j\not=\theta_{i}}\frac{1}{\theta_{i}-y_j}.\ee
It can be shown that for all large $N$, $\rho_N(\theta_{i})\in [L_i(N); D_i(N)]$ and $\rho_N(\theta_{i})=\frac{L_i(N)+ D_i(N)}{2}+o(\frac{1}{\sqrt{N}}).$

To define the limiting correlation function at an outlier, we consider for $k=1, 2, \ldots,$ the distribution $G_k(\cdot )$ given by
\begin{equation}\label{eq:Gkdef}
  G_k(x) = \frac1{Z_k}  \int_{-\infty}^x\cdots \int_{-\infty}^x
  \prod_{1\le i<j\le k} |\xi_i-\xi_j|^2
  \cdot \prod_{i=1}^k e^{-\frac12 \xi_i^2}
  d\xi_1\cdots d\xi_k .
\end{equation}

In other words, $G_k$ is the distribution of the \emph{largest
eigenvalue} of $k\times k$ GUE. 
It has been shown (see \cite{Mehta} or \cite{AGZ09} e.g.) that \begin{equation}\label{eq:Hop}
  G_k(x) = \det(1-\mathbf{H}_x^{(k)}),
\end{equation}
where $\mathbf{H}^{(k)}_x$ is the operator acting on
$L^2((x,\infty))$ defined by the Christoffel Darboux kernel of some rescaled Hermite polynomials 
 satisfying the orthogonality relationship
$\int_{-\infty}^\infty p_m(x)p_n(x) e^{-\frac12 x^2} dx =
  \delta_{mn}$ . We refer the reader to \cite{BBP}, Section 1.2.2 for a more complete statement of this fact.

Let us denote by $\lambda_{max}$ the largest of the $k_i$ outliers around $\rho_N(\theta_{i})$. 
\bt \label{theo: S}There exists $c>0$ depending on $\theta_{i}$ and $\nu$ only such that  $$\lim_{N \to \infty } \mathbb{P}(\sqrt{N}c(\lambda_{max}- \rho_N(\theta_{i})\leq x)=G_{k_i}(x).$$
\et 
We actually prove that the $k_i$ outliers around $\rho_N(\theta_{i})$ fluctuate as the eigenvalues of a $k_i\times k_i$ GUE.
\paragraph{}
Finally, we turn to the fluctuations in a neighborhood of an isolated point of vanishing density.
Let $u_0 \in \mathbb{R}$ be such that $p (u_0)=0$ and that there exists $\epsilon >0$ such that, $\forall u \in ]u_0-\epsilon;u_0+ \epsilon[ \setminus \{u_0\}$,   $p(u) >0$ . 
Assume  that  for any $\theta_i $ such that $\int \frac{d\nu(s)}{(\theta_i -s)^2}=1$, we have $\theta_i +m_\nu(\theta_i) \neq u_0.$
We  prove in Proposition \ref{approxPearcey}, that  for $\eta$ small enough, for all large N,
there exists $u_N $ in  $
]u_0-\eta;u_0+ \eta[$ such that $p_{N}(u_N)=0$ and  $\forall u \in ]u_0-\eta;u_0+ \eta[\setminus \{u_N\}$, $p_{N}(u)>0$, where $p_N$ denotes the density of $\mu_{sc}\boxplus \mu_{A_N}$.\\
Last we derive the asymptotic behavior of eigenvalues at the vicinity of $u_N$.
Consider the
Pearcey kernel defined by 
\be \label{defKP}K_P(x,y):=\frac{1}{2i\pi} \int_{\Gamma_0} dt \int_{-i\infty}^{i\infty}ds e^{-t^4+xt +s^4-sy}\frac{1}{s-t} .\ee
The contour $\Gamma_0$ is formed by two curves lying respectively to the right and left of $0$: one goes from $\infty e^{i\frac{\pi}{4}}$ to $\infty e^{-i\frac{\pi}{4}}$ and the other from $-\infty e^{i \pm \frac{\pi}{4}}$ to $-\infty e^{-i\frac{\pi}{4}}.$ See Figure \ref{fig: contP} below.
\begin{figure}[h!]\begin{center}
\begin{picture}(0,0)%
\includegraphics{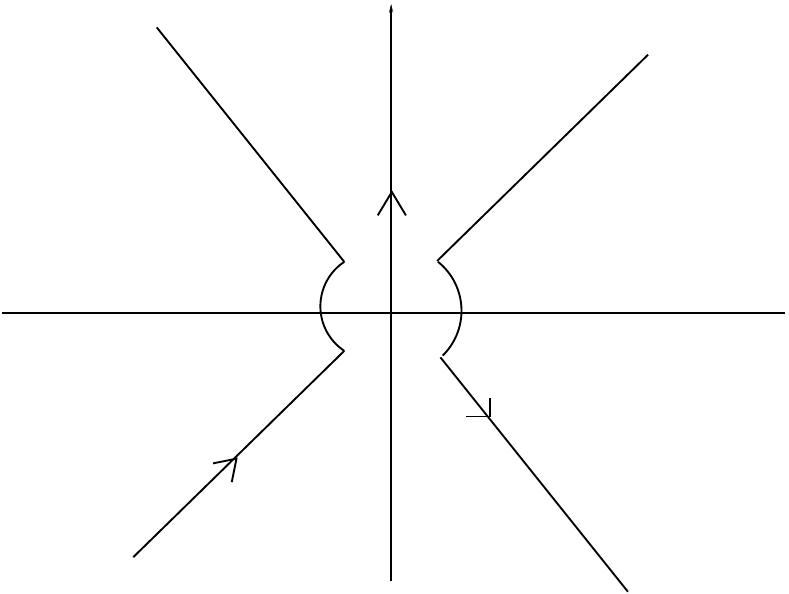}%
\end{picture}%
\setlength{\unitlength}{1184sp}%
\begingroup\makeatletter\ifx\SetFigFont\undefined%
\gdef\SetFigFont#1#2#3#4#5{%
  \reset@font\fontsize{#1}{#2pt}%
  \fontfamily{#3}\fontseries{#4}\fontshape{#5}%
  \selectfont}%
\fi\endgroup%
\begin{picture}(12591,9470)(-32,-8973)
\put(6676,-1561){\makebox(0,0)[lb]{\smash{{\SetFigFont{7}{8.4}{\rmdefault}{\mddefault}{\updefault}{\color[rgb]{0,0,0}$\gamma_0$}%
}}}}
\put(2551,-1561){\makebox(0,0)[lb]{\smash{{\SetFigFont{7}{8.4}{\rmdefault}{\mddefault}{\updefault}{\color[rgb]{0,0,0}$\Gamma_0$}%
}}}}
\end{picture}%
\caption{Contours defining the Pearcey kernel}\label{fig: contP}
\end{center}
\end{figure}
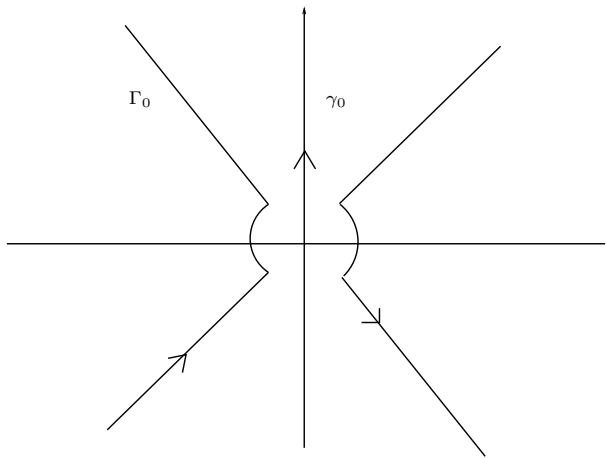
The Pearcey distribution has been defined in \cite{TW2}, \cite{Adleretal1}, \cite{Adleretal2}.
Let $k$ be a fixed integer and $f: \R^k \to R$ be a symmetric bounded function with compact support.
\bt \label{theo: P}There exists $\kappa>0$ such that 
\begin{eqnarray*}&&\E \sum_{\small{1\leq i_1 <i_2< \cdots <i_k \leq N}} f\left (\kappa N^{\frac{3}{4}}(\lambda_{i_1}-u_N),\kappa N^{\frac{3}{4}}( \lambda_{i_2}-u_N), \ldots, \kappa N^{\frac{3}{4}}(\lambda_{i_k}-u_N) \right)\cr
&&\underset{N \to \infty}{\to}\int_{\R^k}\frac{1}{k!} f(x_1, \ldots, x_k) \det(K_P(x_i,x_j))_{i,j=1}^k \prod_{i=1}^k dx_i.
\end{eqnarray*}
\et

The article is organized as follows. In Section \ref{sec: freeconvgauss}, we review the fundamental properties of the free convolution that we later need in the proof. Section \ref{sec:compsupp} gives fine estimates on the comparison of the support of the spectral distribution of $M_N$ on the one hand and that of $\mu_{sc}\boxplus \nu$ on the other hand.
These are the fundamental tools for the asymptotic analysis of eigenvalue statistics in Section \ref{Sec: asanal}. Therein the basic tool is a saddle point analysis of the correlation functions of the deformed GUE.
\section{Free convolution by a semicircular distribution \label{sec: freeconvgauss}}
\subsection{The free convolution}
We recall here an  analytic definition of the free convolution of two probability measures.
Let $\tau $ be a probability measure on $\R$. 
Its Stieltjes transform $m_\tau$ is defined by    
$$m_\tau(z):=\int \frac{1}{z-y}d\tau(y).$$ $m_{\tau}$ is analytic on the complex upper half-plane $\C^+.$
There exists a domain $$D_{\alpha , \beta } = \{ u+iv \in \C, |u| < \alpha v, v > \beta \}$$
on which $m_\tau $ is univalent. 
Let $K_\tau $ be its inverse function, defined on $m_\tau (D_{\alpha , \beta })$, and 
$$R_\tau (z) = K_\tau (z) - \frac{1}{z}.$$
\bdefi Given two probability measures $\tau $ and $\nu $, 
there exists a unique probability measure $\lambda $ such that
$$R_\lambda = R_\tau + R_\nu $$
on a domain where these functions are defined. 
The probability measure $\lambda $ is called 
the free convolution of $\tau $ and $\nu $ and denoted by $\tau \boxplus \nu $.\edefi We refer the reader to \cite{V,VDN, HP, NS} for an introduction to free probability theory.
The free convolution of probability measures has an important property, 
called subordination, which can be stated as follows: 
let $\tau $ and $\nu $ be two probability measures on $\R$; 
there exists an analytic map $\omega: \C^+ \rightarrow \C^+$
 such that $\omega(z)/z \rightarrow 1$ as $z\rightarrow \infty$ with $z\in D_{\alpha,\beta}$, for every such domain, and such that$$\forall z \in \C^+ ,  ~~~~m_{\tau \boxplus \nu}(z)= m_\nu (\omega(z)).$$
This phenomenon was first observed by D. Voiculescu under a genericity assumption in \cite{Voiculescu93}, 
and then proved in generality in \cite{Biane98} Theorem 3.1. 
Later, a new proof of this result was given in \cite{BelBer07}, 
using a fixed point theorem for analytic self-maps of the upper half-plane. \\
In \cite{Biane97b}, P. Biane provides a deep study of the free convolution by a semicircular distribution, based on this subordination property. 
\subsection{The free convolution $\mu _{sc} \boxplus \nu $}

We first recall here some of Biane's results that will be useful in this paper.
\noindent Let $\nu $ be a probability measure on $\R$. 
P. Biane \cite{Biane97b} introduces the set 
$$\Omega _{ \nu }:=\{ u+iv \in \C^+, v > v_{ \nu }(u)\},$$ 
where the function $v_{ \nu }: \R \rightarrow \R^+$ is defined by 
$$v_{ \nu }(u) = \inf \left\{v \geq 0, \int_{\R} \frac{d\nu (x)}{(u-x)^2+v^2} \leq 1\right\},$$ 
and proves the following

\begin{proposition}\cite{Biane97b}\label{homeo} 
The map 
$$H_{ \nu }: z \longmapsto z+m_\nu (z)$$ 
is a homeomorphism from $\overline{\Omega _{ \nu }}$ to $\C^+ \cup \R$ 
which is conformal from $\Omega _{ \nu }$ onto $\C^+$. 
Let $\omega_{ \nu }: \C^+ \cup \R \rightarrow \overline{\Omega _{\nu }}$ 
be the inverse function of $H_{ \nu }$. 
One has, 
$$\forall z \in \C^+,  ~~~~m_{\mu _{sc} \boxplus \nu }(z)= m_{\nu }(\omega_{ \nu }(z))$$
and then 
\begin{equation}\label{Fulton98} 
\omega_{ \nu }(z)=z-m_{\mu _{sc} \boxplus \nu }(z).
\end{equation}
\end{proposition}

\noindent 

~

\noindent 
The previous results of \cite{Biane97b} allows to conclude that 
$\mu _{sc} \boxplus \nu $ is absolutely continuous with respect to the Lebesgue measure 
and to obtain the following description of the support.

\begin{theoreme} \cite{Biane97b} \label{theoBiane}
Define $\Psi _{ \nu }: \R \rightarrow \R$ by: 
$$\Psi _{ \nu }(t)=H_{ \nu }(t+iv_{ \nu }(t)) = 
t+ \int_{\R} \frac{(t-x)d\nu (x)}{(t-x)^2+v_{\nu} (t)^2}.$$
$\Psi _{ \nu }$ is a homeomorphism and, at the point $\Psi _{ \nu }(t)$, 
the measure $\mu _{sc} \boxplus \nu $ has a density given by
\begin{equation}\label{densite} p_{ \nu }(\Psi _{\nu }(t))=\frac{v_{ \nu }(t)}{\pi }.\end{equation}
Define the set 
\begin{equation}\label{defU} U_{ \nu }:= \left\{ t \in \R, \int_{\R} \frac{d\nu (x)}{(t-x)^2} > 1 \right\} = 
\left\{ t \in \R , v_{ \nu }(t) > 0 \right\}.\end{equation}
The support of the measure $\mu _{sc} \boxplus \nu $ is 
the image of the closure of the open set $U_{ \nu }$ by the homeomorphism $\Psi _{ \nu }$. 
$\Psi _{ \nu }$ is strictly increasing on $U_{ \nu }$. 
\end{theoreme}

\noindent 
Hence, 
$$\R\setminus {\rm supp}(\mu _{sc}\boxplus \nu )= 
\Psi _{ \nu } (\R\setminus \overline{U_{\nu }}).$$
One has $\Psi _{ \nu } = H_{\nu }$ 
on $\R\setminus \overline{U_{ \nu }}$ and 
$\Psi ^{-1}_{ \nu }=\omega_{ \nu }$ on $\R\setminus {\rm supp} (\mu _{sc}\boxplus \nu )$. 
In particular, we have the following description of the complement of the support: 
\begin{equation}\label{ComplSupp} 
\R\setminus {\rm supp}(\mu _{sc} \boxplus \nu ) = H_{ \nu }(\R\setminus \overline{U_{ \nu }}).
\end{equation}

\noindent The following result will be useful later on.

\begin{lemme}\cite{Biane97b} \label{lemmeBiane}
If $t_0$ is a point in the complement of the support of $\nu$ where two components of the set $U_{\nu}$ merge into one, then
$$\int \frac{d\nu(x)}{(t_0-x)^2}=1,$$ $$\int \frac{d\nu(x)}{(t_0-x)^3}=0.$$ 

\end{lemme}

In \cite{CDFF}, when $\nu $ is  a compactly supported probability measure,  the authors establish the following results.

\begin{proposition}\label{compconnexes}\cite{CDFF}
\begin{equation}\label{caract} 
\overline{U_{ \nu }} = {\rm supp}(\nu ) \cup 
\{ t \in \R\setminus {\rm supp}(\nu ), \int_{\R} \frac{d\nu (x)}{(t-x)^2} 
\geq 1 \}.
\end{equation}
Each connected component of $\overline{U_{ \nu }}$ 
contains at least one connected component of ${\rm supp}(\nu )$.\\
\end{proposition}

 We also need the following additional basic results.

\begin{lemme}\label{croissancePsi}
Let $]a;b[ \subset  \R \setminus \{ U_\nu \cup \text{supp}(\nu)\}$. Then,  $\Psi_\nu$ is strictly increasing on $]a;b[$.
\end{lemme}
\noindent {\bf Proof:} 
Since $\forall t \in \R \setminus U_\nu,~~ v_\nu(t)=0$, we have $\Psi_\nu= H_\nu$ on $]a;b[$. Moreover  $\forall t \in ]a;b[, ~~H_\nu^{'} (t) =1-\int  \frac{d\nu(x)}{(t-x)^2} \geq 0.$
The result readily follows since moreover $\Psi_\nu$ is one to one. $\Box$

\begin{lemme}\label{prelim}
If $t \notin {\rm supp}(\nu ) $ is such that there exists $\delta>0$ such that 
\begin{equation} \label{t} ]t-\delta;t[ \subset {U_\nu}
\mbox{~~and~~} [t; t+\delta[ \subset  \R \setminus {U_\nu}.\end{equation}
Then, one has that
$$(i): \: \int \frac{d\nu(x)}{(t-x)^2}= 1, \text{ and }
(ii): \: \int \frac{d\nu(x)}{(t-x)^3}>0.$$
If $t'\notin {\rm supp}(\nu ) $ is such that there exists $\delta>0$ such that \begin{equation} \label{s} ]t'-\delta;t'] \subset \R  \setminus {U_\nu}
\mbox{~~and~~} ]t'; t'+\delta[ \subset  {U_\nu}.\end{equation}
Then, one has that $$ (iii): \: \int \frac{d\nu(x)}{(t'-x)^2}= 1\text{ and }
(iv): \: \int \frac{d\nu(x)}{(t'-x)^3}<0.$$
\end{lemme}
\noindent{\bf Proof:} Since $t$ and $t'$ are in $\overline{U_\nu}\setminus {U_\nu}$, 
(\ref{caract}) readily implies (i) and (iii).
Let us establish (ii). Let $\epsilon>0$ be such that $]t-\epsilon; t+\epsilon[ \subset  \R \setminus {\rm supp}(\nu )$. Set $$f:\begin{array}{cc}]t-\epsilon; t+\epsilon[
\rightarrow \mathbb{R}\\
s \mapsto \int \frac{d\nu(x)}{(s-x)^2}. \end{array}$$
Note that $f^{''}(s) =6\int \frac{d\nu(x)}{(s-x)^4}>0$ so that $f^{'}$ is strictly increasing on $]t-\epsilon; t+\epsilon[$.
Therefore if  $ -f^{'}(t) =2 \int \frac{d\nu(x)}{(t-x)^3} \leq 0$ then $f^{'}>0$ on $]t; t+\epsilon[$ and $\int \frac{d\nu(x)}{(s-x)^2}> 1$
for $s\in  ]t; t+\epsilon[$ which leads to a contradiction with (\ref{t}). Similarly, one can prove (iv).$\Box$

\brem
In the rest of the article, since we deal with a measure $\nu$ satisfying (\ref{suppnu}), we have $ {\rm supp}(\nu ) \subset U_\nu$.
\erem

\subsection{\label{subsec:22}The free convolution $\mu _{sc} \boxplus \mu _{A_N}$ and the localization of the spectrum of $M_N$}
In \cite{CDFF}, the authors prove that a precise localization of the spectrum of $M_N$ can be described thanks to the support of the free convolution $\mu _{sc} \boxplus \mu _{A_N}$. In this section, we recall some of their results that we need afterwards.

\begin{theoreme}\label{inclusion}\cite{CDFF} One has that
$\forall \epsilon > 0$, 
$$\mathbb{P}(\text{ For all large N }, {\rm Spect}(M_N)\subset 
\{ x, {\rm dist}(x, {\rm supp}(\mu _{sc} \boxplus \mu _{A_N}))\leq \epsilon \} )=1.$$
\end{theoreme}
An outlier in the spectrum of $M_N$ is an eigenvalue of $M_N$ lying outside the support of $\mu _{sc} \boxplus \nu.$ As we now explain,  
it is possible to describe outliers thanks to the support of $\mu _{sc} \boxplus \mu_{A_N}$.  

\paragraph{Notations and definitions }Throughout the rest of the article, we denote $U_\nu$, $H_\nu$, $\Psi_\nu$, $v_\nu$ and $p_\nu$ by $U$, $H$, $\Psi$, $v$ and $p$ respectively.
We also denote $U_{\mu_{A_N}}$, $H_{\mu_{A_N}}$, $\Psi_{\mu_{A_N}}$, $v_{\mu_{A_N}}$, $p_{\mu_{A_N}}$ by $U_N$, $H_N$, $\Psi_N$, $v_N$ and $p_N$ respectively.  Last, we define the probabily measure  $\hat{\nu}_N$ by $$\hat{\nu}_N= \frac{1}{N-r} \sum_{i=1}^{N-r} \delta_{\beta_i(N)}.$$ It is easy to see that $\hat{\nu}_N$ weakly converges to $\nu$.
We define $$
 \Theta _{ \nu }
=\Theta \cap (\R\setminus \overline{U}).
$$
Furthermore, for any $\theta _j \in \Theta _{ \nu }$, we set 
\begin{eqnarray}
\rho _{\theta _j}:=H(\theta _j)=\theta _j+m_\nu (\theta _j).
\end{eqnarray}
Note that $\rho_{\theta_j}$ lies outside of the support of $\mu _{sc} \boxplus \nu $ according to \eqref{ComplSupp}.
Define also
\begin{eqnarray}
K_{ \nu }(\theta _1, \ldots , \theta _J):=
{\rm supp}(\mu _{sc} \boxplus \nu )\bigcup \left\{ \rho _{\theta _j}, \, \theta _j \in \, \Theta _{ \nu }\right\} .
\end{eqnarray}
In \cite{CDFF}, the authors obtain moreover the following inclusion of the support of $\mu _{sc} \boxplus \mu _{A_N}$.
\begin{theoreme}\label{Theo-InclSupportN}
For any $\epsilon > 0$, 
$${\rm supp}(\mu _{sc} \boxplus \mu _{A_N})\subset 
K_{ \nu }(\theta _1, \ldots , \theta _J) + (-\epsilon, \epsilon),$$ 
when $N$ is large enough.
\end{theoreme}


In \cite{CDFF}, the authors proved this theorem when the $\text{supp} (\nu)$ has a finite number of connected components.
Nevertheless, it is still true in our more general setting as we prove in the following lines.
We will use the following lemma  in \cite{CDFF} which proof does not care about the number of connected components of the supports.
\begin{lemme}{\label{Lemme-InclU-N}}\cite{CDFF}
 For any $\epsilon>0$, 
\begin{eqnarray}{\label{InclusionUN}}
U_{N}\subset  \{ x, \text{dist}(x, \overline{U }) < \epsilon \} 
\cup \, \{ x, \, {\rm dist}(x, \Theta _{ \nu }) < \epsilon \},
\end{eqnarray}
for all large $N$.
\end{lemme}
\noindent {\bf Proof of Theorem \ref{Theo-InclSupportN}:}\\
First,  one can readily observe that if $x$ satisfies 
${\rm dist}(x, {\rm supp}(\nu ))\geq 1 $ 
then  
$-m_{\nu}'(x)\leq 1$. 
This implies that the open set $U$ is included in the compact set 
$\{ x, \, {\rm dist}(x, {\rm supp}(\nu ))\leq 1\}.$
Then we can choose  $K$  large enough such that 
$\{x, \text{dist}(x,\overline{U}\cup \Theta _{ \nu }) \leq 1\} \subset [-K;K]$ and, since $\lim_{y \rightarrow \pm \infty} \Psi(y)=\pm\infty$ and $(\text{supp} (\mu_{A_N} \boxplus \mu_{sc}))_{N}$ are uniformly bounded, \begin{equation}\label{bornes}\text{supp} (\mu_{A_N} \boxplus \mu_{sc}) \subset [\Psi(-(K-1); \Psi(K-1)].\end{equation}
Let $\epsilon>0$. Since $\Psi $ is uniformly continuous on $[-K;K]$, there exists $0<\alpha <1$ such that 
\begin{equation}\label{premsup}
\Psi(\{x, {\rm dist}(x, \overline{U} \cup \Theta_{ \nu } ) < \alpha \}) \subset \{y, \text{dist}(y,K_{ \nu }(\theta _1, \ldots , \theta _J) ) < \epsilon\}.
\end{equation}
Since according to Lemma 
\ref{Lemme-InclU-N}, $\overline{U_{N}}\subset  \{ x, \text{dist}(x, \overline{U }\cup \Theta _{ \nu }) < \alpha/2 \}$ for all large $N$, we have
\begin{equation}\label{secsupp} \Psi_N([-K;K]\cap \{x, \text{dist}(x,  \overline{U }\cup \Theta _{ \nu })\geq \frac{\alpha}{2} \})\subset
\Psi_N(\R \setminus \overline{U_N }) = \R \setminus \text{supp} (\mu_{A_N} \boxplus \mu_{sc}).\end{equation}
Denote by ${\cal A}_{\alpha,K}$ the set $[-K;K] \cap \{x, \text{dist} (x,  \overline{U }\cup \Theta _{ \nu }) \geq \alpha \}.$
Note that for all large $N$, ${\cal A}_{\alpha/2,K} \subset \R \setminus \overline{U_N}$.
Moreover, $$\bigcup_{x\in {\cal A}_{\alpha,K-1} }]x -\frac{\alpha}{2}; x +  \frac{\alpha}{2 }[ \subset   {\cal A}_{\alpha/2,K}.$$
Thus, using Lemma \ref{croissancePsi} for $\Psi_N$, we get that 
$$\bigcup_{x\in {\cal A}_{\alpha,K-1} }]\Psi_N(x -\frac{\alpha}{2}); \Psi_N(x +  \frac{\alpha}{2 })[ \subset  \Psi_N( {\cal A}_{\alpha/2,K}).$$
Now, using the assumptions ($H_3$) on the spectrum of $A_N$, it is easy to see that $\Psi_N$ converges uniformly towards $\Psi$ on the compact set ${\cal A}_{\alpha/2,K}$. 
Moreover, since $\Psi$ is continuous on the compact set ${\cal A}_{\alpha,K-1}$, we have $$\inf_{x \in {\cal A}_{\alpha,K-1}} \min \left( \vert \Psi(x-\alpha/2)-\Psi(x)\vert; \vert \Psi(x+\alpha/2)-\Psi(x)\vert\right) =m>0.$$
Therefore since  for all large $N$, $\sup_{ {\cal A}_{\alpha/2,K}} \vert \Psi_N(x)-\Psi(x) \vert <m$
and  using also Lemma \ref{croissancePsi} for $\Psi$,  we get that for all large $N$, for all $x \in  {\cal A}_{\alpha,K-1}$,
$\Psi_N(x -  \frac{\alpha}{2 })< \Psi(x)<\Psi_N(x +  \frac{\alpha}{2 })$ and therefore
\begin{equation}\label{tersupp} \Psi( {\cal A}_{\alpha,K-1} ) \subset  \Psi_N({\cal A}_{\alpha/2,K}). \end{equation}
(\ref{secsupp}) and (\ref{tersupp}) yield that
$$\text{supp} (\mu_{A_N} \boxplus \mu_{sc}) \subset \R \setminus \Psi(  {\cal A}_{\alpha,K-1})$$
with $$ \R \setminus \Psi(  {\cal A}_{\alpha,K-1})=]-\infty; \Psi(-K+1)[ \cup ]\Psi(K-1); +\infty[
\cup \Psi(  \{x, \text{dist} (x,  \overline{U }\cup \Theta _{ \nu }) < \alpha \}.$$
Then, the result readily follows from (\ref{bornes}) and (\ref{premsup}). $\Box$\\
~~

Moreover, in \cite{CDFF}, the authors proved that the spikes  of the perturbation which belong to $\R\setminus \overline{U}$, generate outliers in the spectrum of the deformed model. 

\begin{theoreme}{\label{ThmASCV}}\cite{CDFF} 
Let  $\theta_j \in \R\setminus \overline{U}$ 
(i.e. $\in \Theta _{ \nu }$).  Denote by $n_{j-1}+1, \ldots , n_{j-1}+k_j$ the descending ranks of $\theta _j$ among the eigenvalues of $A_N$. Then,  almost surely,
$$\lim_{N \to \infty}\lambda_{n_{j-1}+i}(M_N)= \rho _{\theta _j}=H(\theta _j), \: \forall \, 1 \leq i \leq k_j.$$

\end{theoreme}
\brem\label{extension} In \cite{CDFF}, the authors proved this theorem when the support of $\nu$ has a finite number of connected components; nevertheless it is still true in our more general setting since it follows from Theorems \ref{Theo-InclSupportN} and  
\ref{inclusion}  and an exact separation phenomenon (see Theorem 7.1 in \cite{CDFF}) which proof does not care about the number of connected components of the support of $\nu$.
\erem

\section{Comparison of the supports of $\mu_{sc}\boxplus \nu$ and $\mu_{sc}\boxplus\mu_{A_N}$  \label{sec:compsupp}}
As we show in the next Section \ref{Sec: asanal}, the support of  $\mu_{sc}\boxplus\mu_{A_N}$ plays a fundamental role in the study of the fluctuations of eigenvalues at the edges of the spectrum. 
 Due to assumptions $(H_2)$ and $(H_3),$
we are able to show that the supports of $\mu_{sc}\boxplus \nu$ and $\mu_{sc}\boxplus\mu_{A_N}$ exhibit very similar features at edges which are distant from outliers as we explain in Subsection \ref{subsec: UN} below. In Subsection \ref{subsec: outliers}, we prove that $\mu_{sc}\boxplus\mu_{A_N}$ has a connected component in the vicinity of each outlier.
Subsections \ref{subsec: UN1} and \ref{subsec: UN2} are devoted to the proof of the propositions stated in Subsection \ref{subsec: UN} and Subsection \ref{subsec: outliers}.

\subsection{Fundamental preliminary results \label{subsec: UN}}
The two following results will be  fundamental for considering asymptotics of the correlation kernel at the edges of  the support of $\mu_{sc}\boxplus \nu$.
\begin{proposition}\label{rightedge}
Assume that  for a sufficiently small $\epsilon >0$, 
$$p(u)>0, ~~\forall u \in ]u_0-\epsilon; u_0[, ~~\text{and}~~p(u)=0, ~~ \forall  u \in [u_0; u_0+\epsilon[.$$
Set $t_0=\Psi^{-1}(u_0).$ Then there exists $\tau>0$ such that $]t_0-\tau; t_0[ \subset U$, $ [t_0; t_0 +\tau[ \subset \mathbb{R} \setminus U$ and we have  $\int \frac{d\nu(x)}{(t_0 -x)^2}=1$ and $\int \frac{d\nu(x)}{(t_0 -x)^3}>0$. Assume that for all $j\in\{1,\ldots,J\}, \theta_j \neq t_0$.
Then for $\tau>0$ small enough, for all large $N$, there exists one and only one $t_0(N)$ in $]t_0-\tau; t_0+\tau[$, such that $\int \frac{1}{(t_0(N) -x)^2} d\mu_{A_N}(x)=1$ and $]t_0-\tau; t_0+\tau[\cap U_N= ]t_0-\tau; t_0(N)[$. Moreover, for $\eta>0$ small enough, for all large $N$, 
$u_0(N)= \Psi_N(t_0(N))\in ]u_0-\eta;  u_0+\eta[$, $$\text{and}~~\forall u  \in 
  ]u_0-\eta;u_0(N) [ , ~~p_N(u)>0 ~~\text{and} ~~\forall u  \in [u_0(N); u_0+ \eta[, ~~p_N(u)=0.$$

 Moreover, we have \begin{equation}\label{unuo} u_0(N)=u_0  + \epsilon_N(t_0(N)) +\frac{1}{4}({\epsilon_N}^{'}(t_0(N)))^2(1+o(1)) + O(\frac{1}{N}),\end{equation}
 where for $t$ in a small neighborhood of $t_0$,
$$\epsilon_N(t)= \frac{N-r}{N}\int \frac{d\hat{\nu}_N(x)}{(t-x)}-\int \frac{d\nu(x)}{(t-x)}.$$

\end{proposition}
Similarly we have the following result involving the left edges of the support of $\mu_{sc}\boxplus \nu$.
\begin{proposition}\label{leftedge}
Assume that  for a sufficiently small $\epsilon >0$, 
$$p(u)>0, ~~\forall u \in ]u_0; u_0+ \epsilon[, ~~\text{and}~~p(u)=0, ~~ \forall u \in [u_0-\epsilon; u_0[.$$
Set $t_0=\Psi^{-1}(u_0).$    Then there exists $\tau>0$ such that $]t_0-\tau; t_0] \subset \mathbb{R} \setminus U$, $ ]t_0; t_0 +\tau[ \subset  U$ and we have 
 $\int \frac{d\nu(x)}{(t_0 -x)^2}=1$ and $\int \frac{d\nu(x)}{(t_0 -x)^3}<0$. Assume that for all $j\in\{1,\ldots,J\}, \theta_j \neq t_0$.
Then for $\tau>0$ small enough, for all large $N$, there exists one and only one $t_0(N)$ in $]t_0-\tau; t_0+\tau[$, such that $\int \frac{1}{(t_0(N) -x)^2} d\mu_{A_N}(x)=1$ and $]t_0-\tau; t_0+\tau[\cap U_N= ] t_0(N);t_0+\tau[$. Moreover, for $\eta>0$ small enough, for all large $N$, 
$u_0(N)= \Psi_N(t_0(N))
\in ]u_0-\eta;  u_0+\eta[$ $$\text{and}~~\forall u  \in 
  ]u_0(N); u_0+\eta [ , ~~p_N(u)>0 ~~\text{and} ~~\forall u  \in ]u_0-\eta; u_0(N)], ~~p_N(u)=0.$$
 Moreover we have $$u_0(N)=u_0 + \epsilon_N(t_0(N)) +\frac{1}{4}({\epsilon_N}^{'}(t_0(N)))^2 (1+o(1)) + O(\frac{1}{N}),$$
 where for $t$ in a small neighborhood of $t_0$,
$$\epsilon_N(t)= \frac{N-r}{N}\int \frac{d\hat{\nu}_N(x)}{(t-x)}-\int \frac{d\nu(x)}{(t-x)}.$$

\end{proposition}

\brem
It is clear that, under the assumption (\ref{Scherbina}) of Shcherbina  (\cite{Shcherbina2}),  Theorem \ref{theo: A} and (\ref{unuo}) imply her result.
\erem

The following proposition will be fundamental to study the asymptotics of the correlation kernel in a neighborhood of any  point of
the support of $\mu_{sc}\boxplus \nu$ where the density vanishes.
\begin{proposition}\label{approxPearcey}
Let $u_0 \in \mathbb{R}$ be such that $p (u_0)=0$ and  there exists $\epsilon >0$ such that, $\forall  u \in ]u_0-\epsilon;u_0+ \epsilon[ \setminus \{u_0\}$,   $p(u) >0$ . 
\\Set $t_0=\Psi^{-1} (u_0) \in \R$. Then $t_0$ is a point in $\R \setminus \text{supp}(\nu)$ where two components of $U$ merge and satisfies
$\int \frac{d\nu(s)}{(t_0-s)^2}=1$,  $ \int \frac{d\nu(s)}{(t_0-s)^3}=0.$ We have $u_0=H(t_0)$.\\
Assume  that  for any $i=1,\dots,J$, $\theta_i \neq t_0$. 
Then, for $\eta$ small enough, for all large N,
there exists $u_0(N) $ in  $
]u_0-\eta;u_0+ \eta[$ such that $p_{N}(u_0(N))=0$ and  $\forall u \in ]u_0-\eta;u_0+ \eta[\setminus \{u_0(N)\}$, $p_{N}(u)>0$.\\
  $t_0(N)=\Psi_N^{-1} (u_0(N))$ is a point of $\R \setminus \text{Spect}(A_N) $ where two components of $U_N$ merge and satisfies $\int \frac{d{\mu_{A_N}}(s)}{(t_0(N)-s)^2}=1$,  $\int \frac{d{\mu_{A_N}}(s)}{(t_0(N)-s)^3}=0$. We have  $u_0(N)= H_N(t_0(N)) $
and  $\lim_{N \rightarrow + \infty}t_0(N)=t_0. $

\end{proposition}

\subsection{In the vicinity of outliers \label{subsec: outliers}}
It turns out that the support of $\mu_{sc}\boxplus\mu_{A_N}$ exhibits a small connected component in the vicinity of each outlier.
\begin{proposition}\label{outlier} Let $\theta_i$ be such that $\int \frac{d\nu(x)}{(\theta_i-x)^2} <1$ and $\rho_{\theta_i}=H(\theta_i)$. Then, for $\epsilon>0$ small enough, for all large $N$, $\text{supp}(\mu_{sc}\boxplus\mu_{A_N})$ has a unique connected component $[L_i(N); D_i(N)]$ inside $]\rho_{\theta_i} -\epsilon;  \rho_{\theta_i} +\epsilon[$.
Moreover, setting $\rho_N(\theta_i)= \frac{1}{N} \sum_{y_j \neq \theta_i} \frac{1}{\theta_i -y_j} +\theta_i$, we have 
$$L_i(N)=\rho_N(\theta_i) -2 \sqrt{k_i} \sqrt{1-\int \frac{1}{(\theta_i -x)^2}d\nu(x)}\frac{1}{\sqrt{N}} +o(\frac{1}{\sqrt{N}}),$$
$$D_i(N)=\rho_N(\theta_i) +2 \sqrt{k_i} \sqrt{1-\int \frac{1}{(\theta_i -x)^2}d\nu(x)}\frac{1}{\sqrt{N}} +o(\frac{1}{\sqrt{N}}).$$
Thus, $\rho_N(\theta_i)=\frac{L_i(N)+D_i(N)}{2} +o(\frac{1}{\sqrt{N}}).$
\end{proposition}
\subsection{Some technical lemmas \label{subsec: UN1}}

In the proof of the previous propositions, we will use the following lemmas.

\begin{lemme}\label{montel}
Let $ [a;b] \subset \R \setminus \text{supp}(\nu)\cup \Theta.$ Then 
$m_N: z \mapsto \int \frac{d{\mu_{A_N}}(s)}{(z-s)}$ (resp. $-m_N^{'}: z \mapsto \int \frac{d{\mu_{A_N}}(s)}{(z-s)^2}$) converges  uniformly towards $m: z \mapsto \int \frac{d{\nu}(s)}{(z-s)}$ (resp. $-m^{'}: z \mapsto \int \frac{d{\nu}(s)}{(z-s)^2}$)  on every compact set included in  $\{z \in \mathbb{C}; a < \Re z < b \}$.
\end{lemme}
\noindent {\bf Proof:}
 Let $\gamma>0$ be such that $[a-3\gamma;b+3\gamma] \subset \R \setminus \text{supp}(\nu)\cup \Theta.$
 Since $A_N$ has $N-r$ eigenvalues $\beta _j(N)$ satisfying 
$\max _{1\leq j\leq N-r} {\rm dist}(\beta _j(N),{\rm supp}(\nu ))\vers _{N \rightarrow \infty } 0,$ and the other eigenvalues of $A_N$ are the spikes $\theta_j \in \Theta$,
we can readily deduce that  for all large N, $$[a-2\gamma;b+ 2\gamma] \subset \mathbb{R} \setminus \text{Spect}(A_N) .$$
\noindent It is clear that the  functions $m_N$, $g$, $-m_N^{'}$ and$
-m^{'}$ are holomorphic on $\{z \in \mathbb{C}; a -\gamma< \Re z <b + \gamma \}$.
Since for large $N$,  $\{z \in \mathbb{C}; a-\gamma < \Re z <b +\gamma\} $ is included in $ \{z \in \mathbb{C}; \text{dist}(z;\text{supp}(\nu)) > \gamma;  \text{dist}(z;\text{Spect}(A_N)) > \gamma\}$, it readily follows that for large $N$,  $m_N$ and $m$ (respectively $m_N^{'} $ and $m^{'}$ ) are uniformly bounded by $1/\gamma $ (respectively $1/\gamma^2 $). 
Since the sequence of measures $\mu_{A_N}$ weakly converges to $\nu$, 
 it is easy to see that  $m_N(z)$ (respectively $ m_N^{'}(z)$) converges towards $m(z)$ (respectively $m^{'}(z)$)  for all $z \in  ]a;b[$.
Therefore, by Montel's theorem, the convergence is uniform on every compact set of $\{z \in \mathbb{C}; a -\gamma < \Re z <b +\gamma \} ~~\Box$

\begin{lemme}\label{prelim2}
$(1)$For any $t$ in $U$, $ v_N(t)$ converges towards $v(t)$ when $N$ goes to infinity.\\
$(2)$ For any $t$ in $U$ such that $t \in \R \setminus \{\text{supp}(\nu) \cup \Theta\}$, $\Psi_N(t) $ converges towards $\Psi(t)$ when $N$ goes to infinity.\\
$(3)$  For any $t$ in $ \R \setminus \{\overline{U}\cup \Theta\}$, $\Psi_N(t) $ converges towards $\Psi(t)$ when $N$ goes to infinity.
\end{lemme}
\noindent {\bf Proof :}
Let $t$ be in $U$. Therefore we have $v(t)>0$. Let $0<\epsilon< v(t)$.
We have $\int \frac{d\nu(s)}{(t-s)^2+ (v(t)-\epsilon)^2}>1$ and  $\int \frac{d\nu(s)}{(t-s)^2+ (v(t)+\epsilon)^2}<1$ which implies that for all large $N$, 
$\int \frac{d\mu_{A_N}(s)}{(t-s)^2+ (v(t)-\epsilon)^2}>1$ and  $\int \frac{d\mu_{A_N}(s)}{(t-s)^2+ (v(t)+\epsilon)^2}<1$.
It follows that for all large $N$, $v(t)-\epsilon< v_N(t) < v(t)+\epsilon.$\\
Now, let $t$ be  in $U$  and such that $t \in \R \setminus \{\text{supp}(\nu) \cup \Theta\}$. Let $\delta>0$ such that $[t-\delta; t+ \delta] \subset \R \setminus \{\text{supp}(\nu) \cup \Theta\}.$
 According to Lemma \ref{montel}, $z\mapsto \int \frac{d\mu_{A_N}(x)}{z-x}$ converges towards $z\mapsto \int \frac{d\nu(x)}{z-x}$ uniformly on every compact set of 
$\{t-\delta<\Re z < t+\delta\}.$
Since $v_N(t)$ converges towards $v(t)$, for all large $N$, $ 0 \leq  v_N(t) \leq v(t) +1$. The convergence of $\Psi_N(t) $  towards $\Psi(t)$ when $N$ goes to infinity readily follows from the uniform convergence of 
$z\mapsto \int \frac{d\mu_{A_N}(x)}{z-x}$  towards $z\mapsto \int \frac{d\nu(x)}{z-x}$  on the compact set $\{z =t+i b, 0\leq b \leq v(t)+1 \}$.\\
Let $t$ be in $ \R \setminus \{\overline{U}\cup \Theta\}$. Since $v(t)=0$, we have $\Psi(t)=H(t)$.
Since we assume that $\text{supp}(\nu) \subset U$, we have $t \in \R \setminus \text{supp}(\nu).$ According to the assumption ($H_3$) on the spectrum of $A_N$, for  $\delta >0$ small enough, for all large $N$, we have  $[t-\delta; t+\delta] \subset 
\R \setminus \{ \text{supp}(\nu) \cup \text{supp}(\mu_{A_N})\}$.
Therefore it is easy to see that $\int \frac{d\mu_{A_N}(x)}{(t-x)}$  converges towards $\int \frac{d\nu(x)}{(t-x)}$ and $\int \frac{d\mu_{A_N}(x)}{(t-x)^2}$  converges towards $\int \frac{d\nu(x)}{(t-x)^2} <1$ and thus for all large $N$, $t \notin U_N$.
It follows that for all large $N$, $v_N(t)=0$ and $\Psi_N(t) =H_N(t)=t+\int \frac{d\mu_{A_N}(x)}{(t-x)}$ converges towards $t+\int \frac{d\nu(x)}{(t-x)}=H(t)=\Psi(t)$.
 $\Box$

 \begin{lemme}\label{dk}  Let  $t_0$ $\notin {\rm supp}(\nu ) \cup \Theta$ be  such that $\int \frac{1}{(t_0-x)^2}d\nu(x)=1$ and $\int \frac{1}{(t_0-x)^3}d\nu(x)\neq0$. Then for small enough $\epsilon>0 $, for all large $N$,  there exists one and only one $t_0(N) \in  ]t_0 -\epsilon;t_0 + \epsilon[$ such that 
$\int \frac{1}{(t_0(N) -x)^2} d\mu_{A_N}(x)=1$.
 $t_0(N)$  satisfies
$$t_0(N)=t_0+f_N(t_0(N))$$ where\\

$f_N(t)=$ $$h(t)\left[
 \left\{\frac{N-r}{N}  \int \frac{1}{(t-x)^2}d\hat{\nu}_N(x)- \int \frac{1}{(t-x)^2}d\nu(x)\right\}
 + \frac{1}{N} 
 \sum_{j=1}^J \frac{k_j}{(t-\theta_j)^2}\right]$$
with $h(t) =\frac{1}{\int \frac{(t-x+t_0-x)}{(t-x)^2(t_0-x)^2}d\nu(x)} $ and  $0< K_1(\epsilon) < \vert h(t) \vert< K_2(\epsilon), \forall t \in ]t_0 -\epsilon; t_0+ \epsilon[.$ \\
Moreover,
if $]t_0 -\epsilon;t_0[ \subset U$ and $]t_0; t_0 +\epsilon[ \subset \R \setminus U$ (respectively $
]t_0; t_0+ \epsilon[ \subset U$ and $]t_0-\epsilon; t_0 [ \subset \R \setminus U$) , then for all large $N$, $]t_0 -\epsilon;t_0 + \epsilon[\cap U_N=]t_0 -\epsilon;t_0 (N)[$ (respectively  $]t_0 -\epsilon;t_0 + \epsilon[\cap U_N=]t_0(N);t_0+ \epsilon[$.)
\end{lemme}
\noindent {\bf Proof:}   One can readily see that $t \notin\{\theta_i,i=1,\ldots, J, \beta_j, j=1,\ldots,  N-r\}$ is in $U_N $ if and only if 
 $P_N(t) >0$ where $P_N(t)$ is the polynomial  defined by
 \begin{eqnarray}P_N(t)&=&\prod_{ i=1}^{N-r} (t-\beta_i)^2  \prod_{j=1}^J ( t-\theta_j)^2 \left( \int \frac{d\mu_{A_N}}{(t-x)^2} -1\right)
\\ &=&\frac{1}{N}\sum_{i=1}^{N-r} \prod_{l\neq i} (t-\beta_l)^2  \prod_{j=1}^J (t-\theta_j)^2 \nonumber \\
&&+  \frac{1}{N} \prod_{i=1}^{N-r} (t-\beta_i)^2 \sum_{ j=1^{J}} k_j \prod_{l\neq j}( t-\theta_l)^2 \nonumber\\ && - \prod_{j=1}^J ( t-\theta_j)^2 \prod_{i=1}^{N-r} (u-\beta_i)^2. \label{polynome}\end{eqnarray}
Condition ($H_3$) on the spectrum of $A_N$ allows us to 
choose $\epsilon >0$ small enough  such that  for  $N$ large enough $ [t_0-2\epsilon; t_0 +2\epsilon]$ is in the complement of the support of $\nu$ and the support of $\mu_{A_N}$.
 $ P_N(t)=0$ for $t \in ]t_0-\epsilon;t_0+ \epsilon[$ if and only if 
 \begin{equation}\label{zeros}1 - \frac{N-r}{N}  \int \frac{1}{(t-x)^2}d\hat{\nu}_N - \frac{1}{N} 
 \sum_{j=1}^J \frac{k_j}{( t-\theta_j)^2}=0.\end{equation}
 Using that  
 $$ \int \frac{1}{(t_0-x)^2}d\nu(x)=1,$$
 (\ref{zeros}) can be rewritten as follows:
\begin{eqnarray*} 
&& \int \frac{1}{(t_0-x)^2}d\nu(x)- \int \frac{1}{(t-x)^2}d\nu(x)=\cr
&&\left\{\frac{N-r}{N}  \int \frac{1}{(t-x)^2}d\hat{\nu}_N(x)- \int \frac{1}{(t-x)^2}d\nu(x)\right\}
 + \frac{1}{N} 
 \sum_{j=1}^J \frac{k_j}{( t-\theta_j)^2},\end{eqnarray*}
 or equivalently 

\begin{eqnarray*} 
&&(t-t_0)\int \frac{(t-x+t_0-x)}{(u-x)^2(t_0-x)^2}d\nu(x)=\cr
&&\left\{\frac{N-r}{N}  \int \frac{1}{(t-x)^2}d\hat{\nu}_N(x)- \int \frac{1}{(t-x)^2}d\nu(x)\right\}
 + \frac{1}{N} 
 \sum_{j=1}^J \frac{k_j}{( t-\theta_j)^2}.\end{eqnarray*}
Since 
we have $\int \frac{1}{(t_0-x)^3}d\nu(x)\neq 0$,   it readily follows that for $\epsilon>0$ small enough and for all $z $ such that 
$\vert z-t_0\vert  \leq \epsilon$, 
$\int \frac{(z-x+t_0-x)}{(z-x)^2(t_0-x)^2}d\nu(x)\neq0$. Therefore,  there exists $C_1(\epsilon)>0$ and $C_2(\epsilon)>0$
 such that for any $z$ such that $ \vert z-t_0\vert \leq \epsilon$, $0< C_1(\epsilon)<\vert \int \frac{(z-x+t_0-x)}{(z-x)^2(t_0-x)^2}d\nu(x) \vert< C_2(\epsilon).$
Define on $\{z; \vert z-t_0\vert \leq \epsilon\}$,  $$h(z) =\frac{1}{\int \frac{(z-x+t_0-x)}{(z-x)^2(t_0-x)^2}d\nu(x)}.$$ 
 Using Lemma \ref{montel}, by Rouch\'e theorem, for large $N$, the function $$z-t_0-h(z)\left[
 \left\{\frac{N-r}{N}  \int \frac{d\hat{\nu}_N(x)}{(z-x)^2}- \int \frac{d\nu(x)}{(z-x)^2}\right\}
 + \frac{1}{N} 
 \sum_{j=1}^J \frac{k_j}{\vert z-\theta_j\vert^2}\right]$$ has exactly one zero $z_0$ in $\{z; \vert z-t_0 \vert < \epsilon\}$.
 Since $\bar{z_0}$ is obviously a zero too, we can conclude that $z_0$ is real.
 Hence, for $\epsilon$ small enough,  for all large $N$, $P_N$ has exactly one zero $t_0(N)$ in $]t_0-\epsilon;t_0+\epsilon[$ and 
  \begin{eqnarray*}
 t_0(N)&=&t_0 + h(t_0(N))\left[
 \left\{\frac{N-r}{N}  \int \frac{d\hat{\nu}_N(x)}{(t_0(N)-x)^2}- \int \frac{d\nu(x)}{(t_0(N)-x)^2}\right\}\right. \\ &&\left.
 + \frac{1}{N} 
 \sum_{j=1}^J \frac{k_j}{( t_0(N)-\theta_j)^2}\right]
 \end{eqnarray*}
 where $0< K_1(\epsilon) < \vert h(t_0(N)) \vert < K_2(\epsilon).$

 Now, if $]t_0 -\epsilon;t_0[ \subset U$ and $]t_0; t_0 +\epsilon[ \subset \R \setminus U$ (respectively $
]t_0; t_0+ \epsilon[ \subset U$ and $]t_0-\epsilon; t_0 [ \subset \R \setminus U$), then since  for all large $N$, $P_N(t_0- \epsilon/2)>0$ and $P_N(t_0 +\epsilon/2) <0$ (respectively $P_N(t_0- \epsilon/2)<0$ and $P_N(t_0 +\epsilon/2) >0$), 
it is clear that for all large $N$, $]t_0 -\epsilon;t_0 + \epsilon[\cap U_N=]t_0 -\epsilon;t_0 (N)[$ (respectively  $]t_0 -\epsilon;t_0 + \epsilon[\cap U_N=]t_0(N);t_0+ \epsilon[$.)
The proof of Lemma \ref{dk}  is complete. $\Box$

\begin{lemme}\label{prelim22}
Let $x_0$ be such that $\int \frac{d\nu(s)}{(x_0-s)^2}=1$, $x_0 \neq \theta_j, \forall 1\leq j \leq J$, and  there exists $\tau >0$ such that, $\forall x \in ]x_0-\tau;x_0+ \tau[ \setminus \{x_0\}$,   $\int \frac{d\nu(s)}{(x-s)^2}>1$ .
Then, $ x_0 \notin \text{supp}(\nu) \cup \Theta.$ Set $d_1 =\sup\{s \in \text{supp}(\nu) \cup \Theta; s < x_0\}$ and $d_2=\inf\{s \in \text{supp}(\nu) \cup \Theta; s > x_0\}$.
Let  $[a;b]$ be  such that $x_0 \in ]a;b[$, $[a;b] \subset ]d_1;d_2[.$ Then, $\forall x \in [a;b]\setminus \{x_0\}, ~~\int \frac{d\nu(s)}{(x-s)^2}>1.$
Moreover, for all large N, $[a;b] \subset \mathbb{R} \setminus \text{Spect}(A_N) $ and  there exists $x_0(N) $ in  $
[a;b]$ such that $\int \frac{d{\mu_{A_N}}(s)}{(x_0(N)-s)^2}=1$,  $\int \frac{d{\mu_{A_N}}(s)}{(x_0(N)-s)^3}=0$ and  $\forall x \in [a;b]\setminus \{x_0(N)\}$, $\int \frac{d{\mu_{A_N}}(s)}{(x-s)^2}>1$ .
We have also $\lim_{N \rightarrow + \infty}x_0(N)=x_0. $ 
\end{lemme}

\noindent{\bf Proof}: Since we assume that for any $x$ in $\text{supp}(\nu)$, $\int \frac{d\nu(s)}{(x-s)^2}>1$ (i.e $\text{supp}(\nu) \subset U$) and that  $x_0 \neq \theta_j, \forall 1\leq j \leq J$, it readily  follows that $x_0 \notin  \text{supp}(\nu)\cup \Theta.$ 
Let  $[a;b]$ be  such that $x_0 \in ]a;b[$, $[a;b] \subset ]d_1;d_2[.$
Since $\int \frac{d\nu(s)}{(x_0-s)^2}=1$ and there exists $\tau >0$ such that, $\forall x \in ]x_0-\tau;x_0+ \tau[ \setminus \{x_0\},~~\int \frac{d\nu(s)}{(x-s)^2}>1$,  the strict convexity of $ z \mapsto \int \frac{d{\nu}(s)}{(z-s)^2}$ on $[a;b]$ implies  that   \begin{equation}\label{positive}\forall x \in [a;b] \setminus \{x_0\}, \int \frac{d\nu(s)}{(x-s)^2}>1.\end{equation}
By Lemma \ref{montel}, $\phi_N: z \mapsto \int \frac{d{\mu_{A_N}}(s)}{(z-s)^2}-1$ converges  uniformly towards $\phi: z \mapsto \int \frac{d{\nu}(s)}{(z-s)^2}-1$ on every compact set of $\{z \in \mathbb{C}; a < \Re z <b \} $.\\
By the principle of isolated zeroes, there exist  $\delta_0$ such that   $[x_0 -\delta_0; x_0+\delta_0] \subset ]a;b[$ and $\phi$ has no other zero in $\{z \in \mathbb{C};
\vert z-x_0\vert \leq  \delta_0\}$ than $x_0$. 
Thus,  using Hurwitz's theorem, we can claim that for all large $N$, $\phi_N $ has a unique zero in $\{z \in \mathbb{C};
\vert z-x_0\vert <  \delta_0\}$.  Let us denote this zero by $x_0(N)$. Since $\phi_N(\overline{x_0(N)})=\overline{\phi_N(x_0(N))}=0$, it follows that $x_0(N) $ is a real number.\\
Note that since by (\ref{positive}), $\phi(x_0 -\delta_0)>0$ and $\phi(x_0 +\delta_0) >0$, for all large N we have $\phi_N(x_0 -\delta_0)>0$ and $\phi_N(x_0 +\delta_0) >0$; since $\phi_N $ has a unique zero $x_0(N) $ in $\{z \in \mathbb{C};
\vert z-x_0\vert <  \delta_0\}$, this implies that $\phi_N$ reaches its minimum value in $]x_0 -\delta_0; x_0+\delta_0[$ at $x_0(N)$ and thus $\phi'_N(x_0(N))=0$. Moreover since $\phi_N$ is strictly convex on $[a;b]$, we have $\forall x \in [a;b]\setminus \{x_0(N)\}$, $\phi_N(x)>0$.
Now, similarly, for any $0<\delta< \delta_0$, using Hurwitz's theorem, we can claim that for all large $N$, $\phi_N $ has a unique zero in $\{z \in \mathbb{C};
\vert z-x_0\vert <  \delta\}$ and  for $N \geq N(\delta)$, this unique zero is equal to $x_0(N)$; the convergence of $x_0(N)$ towards $x_0$ follows.  $\Box$

\begin{lemme}\label{spike}  For each i such that $\int \frac{1}{(\theta_i -x)^2}d\nu(x)  < 1$, for $\epsilon>0$ small enough, for all large $N$,  $U_N \bigcap ]\theta_i-\epsilon; \theta_i + \epsilon [= ]t_1^i(N),t_2^i(N)[$ 
where $t_1^i(N)$ and $t_2^i(N)$ satisfy
$$t_1^i(N)= \theta_i  - \sqrt{\frac{k_i}{N} \phi_N(t_1^i(N))}$$
$$t_2^i(N)= \theta_i  + \sqrt{\frac{k_i}{N} \phi_N(t_2^i(N))}$$
with $\phi_N(t)= \frac{1}{1 - \frac{N-r}{N} \int \frac{1}{(t-x)^2}d\hat{\nu}_N (x)- \frac{1}{N} 
 \sum_{j\neq i} \frac{k_j}{( t-\theta_j)^2}}$ 
and $1 \leq \phi_N(t) \leq K(\epsilon)$ for any $t \in ]\theta_i-\epsilon; \theta_i + \epsilon [$ .
\end{lemme}
\noindent {\bf Proof:}
 \noindent Let $\theta_i$ be such that $\int \frac{d\nu(x)}{(\theta_i -x)^2 }< 1$.  Let $\epsilon>0$ be such that $]\theta_i- 4\epsilon; \theta_i+ 4\epsilon[ \subset \R \setminus \left\{ \text{supp}(\nu) \cup\{\theta_j, j \neq i\} \right\}$ and 
$\inf_{z\in \C, \vert z-\theta_i\vert \leq 2\epsilon} \vert \int \frac{d\nu(x)}{(z-x)^2 }- 1 \vert =m \neq 0$. In particular, we have that for any $t$ in $[\theta_i- \epsilon; \theta_i+ \epsilon]$, $\int \frac{d\nu(x)}{(t-x)^2 }< 1$ . According to the assumption ($H_3$) on the spectrum of $A_N$, for all large $N$, $ [\theta_i- 3\epsilon;\theta_i[\cup]\theta_i; \theta_i+3 \epsilon] \subset \R \setminus \text{Spect}(A_N).$
 Note that since $\int \frac{d\mu_{A_N}(x)}{(\theta_i\pm\epsilon-x)^2 }$ converges towards $\int \frac{d\nu(x)}{(\theta_i \pm \epsilon-x)^2 }$, we have moreover  for all large $N$, $\int \frac{d\mu_{A_N}(x)}{(\theta_i \pm \epsilon-x)^2 }<1$, whereas $\int \frac{d\mu_{A_N}(x)}{(\theta_i-x)^2 }=+\infty$.
Therefore, for all large $N$, there exists at least one $s_N \in ]\theta_i- \epsilon; \theta_i[$ and at least one $t_N \in ]\theta_i; \theta_i + \epsilon[$ such that $\int \frac{d\mu_{A_N}(x)}{(s_N-x)^2 }=1$ and $\int \frac{d\mu_{A_N}(x)}{(t_N-x)^2 }=1.$
 Let us study the zeroes of the polynomial $P_N$ defined by (\ref{polynome}) in $\{z; \vert z- \theta_i\vert <\epsilon\}$. We know that there are at least two real zeroes $s_N$ and $t_N$.
 Let us rewrite 
  \begin{eqnarray*}P_N(t)&= &\frac{1}{N}\sum_{i=1}^{N-r} \prod_{l\neq i} (t-\beta_l)^2  \prod_{j=1}^J ( t-\theta_j)^2  +  \frac{1}{N} \prod_{l=1}^{N-r} (t-\beta_l)^2 \sum_{ j\neq i} k_j \prod_{p\neq j}( t-\theta_p)^2 \\&&+  \frac{1}{N} \prod_{j=1}^{N-r} (t-\beta_j)^2  k_i \prod_{l\neq i}( t-\theta_l)^2
 - \prod_{j=1}^J ( t-\theta_j)^2 \prod_{l=1}^{N-r} (t-\beta_l)^2.\end{eqnarray*}
 $ P_N(t)=0$ for $t$ such that  $\vert t-\theta_i\vert <2 \epsilon$ if and only if 
 $$( t-\theta_i )^2 \left\{ 1 - \frac{N-r}{N} \int \frac{1}{(t-x)^2}d\hat{\nu}_N(x) - \frac{1}{N} 
 \sum_{j\neq i} \frac{k_j}{( t-\theta_j)^2} \right\} =  \frac{k_i}{N},$$
Since for all large $N$, $ [\theta_i- 3\epsilon; \theta_i+3 \epsilon] \subset \R \setminus \{\text{supp}(\hat{\nu}_N)\cup \text{supp}(\nu)\} $,  using the same arguments as in the proof of Lemma \ref{montel}, we get easily the uniform convergence on 
any compact set included in $\{z\in \C, \theta_i -3 \epsilon < \Re z < \theta_i+ 3 \epsilon \}$ of 
$ z\mapsto 1 - \frac{N-r}{N} \int \frac{1}{(z-x)^2}d\hat{\nu}_N(x) - \frac{1}{N} 
 \sum_{j\neq i} \frac{k_j}{( z-\theta_j)^2}$   towards $ z \mapsto 1 -  \int \frac{1}{(z-x)^2}d{\nu}(x) $. 

 \noindent Hence, we have for all large $N$,
 $$\inf_{z\in \C, \vert z-\theta_i\vert \leq 2\epsilon} \left|
 1 - \frac{N-r}{N} \int \frac{1}{(z-x)^2}d\hat{\nu}_N(x) - \frac{1}{N} 
 \sum_{j\neq i} \frac{k_j}{( z-\theta_j)^2}\right| \geq m/2$$
and $$\inf_{t \in [\theta_i-2\epsilon; \theta_i+2\epsilon]} \left\{
 1 - \frac{N-r}{N} \int \frac{1}{(t-x)^2}d\hat{\nu}_N (x)- \frac{1}{N} 
 \sum_{j\neq i} \frac{k_j}{( t-\theta_j)^2}\right\} \geq m/2$$
and the zeros of $P_N$  in $\{z\in \C, \vert z-\theta_i\vert <2\epsilon\}$
are the solutions of the equation
$( z- \theta_i )^2 = \frac{k_i}{N} \phi_N(z)$ where $$\phi_N(z)= \frac{1}{1 - \frac{N-r}{N} \int \frac{1}{(z-x)^2}d\hat{\nu}_N (x)- \frac{1}{N} 
 \sum_{j\neq i} \frac{k_j}{( z-\theta_j)^2}}$$ and $0 <\vert \phi_N(z) \vert \leq \frac{2}{m}.$ Therefore, by Hurwitz theorem, for all large $N$, $P_N$ has exactly two zeroes in  $\{z\in \C, \vert z-\theta_i\vert < \epsilon\}$.
Since we have already seen that $P_N$ has at least one zero in $ ]\theta_i- \epsilon; \theta_i[$ and at least one zero in $ ]\theta_i; \theta_i + \epsilon[$, we can conclude that  for all large N, $P_N$ has exactly one  zero $t_1^i(N)$ in  $]\theta_i-\epsilon; \theta_i [$ and  one zero $t_2^i(N)$ in $]\theta_i; \theta_i+\epsilon[$.
Moreover since $\phi_N(t)>0$ on $ [\theta_i-\epsilon; \theta_i+\epsilon]$, we have 
 $$ t_1^i(N)=\theta_i  - \sqrt{\frac{k_i}{N} \phi_N(t)} ~~\text{and}~~  t_2^i(N)= \theta_i + \sqrt{\frac{k_i}{N} \phi_N(t)}.$$

Now, since $P_N(\theta_i)  >0$, it is clear that $U_N \bigcap  ]\theta_i-\epsilon; \theta_i + \epsilon [= ]t_1^i(N),t_2^i(N)[$ .
The proof of  Lemma \ref{spike}  is complete.$\Box$

\subsection{Proof of Propositions  \ref{rightedge}, \ref{leftedge}, \ref{approxPearcey} and \ref{outlier} \label{subsec: UN2}}
\paragraph{ Proof of Proposition \ref{rightedge}:}
Using (\ref{densite}) and (\ref{defU}), it is clear that $\Psi^{-1}(]u_0-\epsilon;u_0[) \subset U$ and  $\Psi^{-1}([u_0;u_0+\epsilon[) \subset \R \setminus U$. 
Note that since we assume that $\text{supp}(\nu) \subset U$, this implies that  $t_0=\Psi^{-1}(u_0) \in \R \setminus \text{supp}(\nu).$ Let $0< \delta< \epsilon$ be such that 
$\Psi^{-1}(]u_0 -\delta; u_0 +\delta[) \subset \R \setminus \{ \text{supp}(\nu)\cup \Theta\}.$ Since according to Theorem \ref{theoBiane}, the homeomorphism $\Psi$ is strictly increasing on $U$, we have 
$\Psi^{-1}(]u_0 -\delta; u_0 [) =]\Psi^{-1}(u_0 -\delta); t_0[\subset U$. Moreover according to Lemma \ref{croissancePsi}, $\Psi^{-1}([u_0; u_0+\delta [) =[t_0; \Psi^{-1}(u_0+\delta)[ \subset \R \setminus U$. 
Thus, according to Lemma \ref{prelim} (i) and (ii), we have  $\int \frac{d\nu(x)}{(t_0-x)^2}= 1$,
$\int \frac{d\nu(x)}{(t_0-x)^3}>0.$ Then, using Lemma \ref{dk}, for $\tau$ small enough, for all large $N$ there exists 
one and only one $t_0(N) \in  ]t_0 -\tau;t_0 + \tau[$ such that 
$\int \frac{1}{(t_0(N) -x)^2} d\mu_{A_N}(x)=1$.
 $t_0(N)$  satisfies
$$t_0(N)=t_0+f_N(t_0(N))$$ where
$$f_N(t)=h(t)\left[
 \left\{\frac{N-r}{N}  \int \frac{d\hat{\nu}_N(x)}{(t-x)^2}- \int \frac{d\nu(x)}{(t-x)^2}\right\}
 + \frac{1}{N} 
 \sum_{j=1}^J \frac{k_j}{(t-\theta_j)^2}\right]$$
with $h(t) =\frac{1}{\int \frac{(t-x+t_0-x)}{(t-x)^2(t_0-x)^2}d\nu(x)} $ and  $0< K_1(\tau) < \vert h(t) \vert< K_2(\tau), \forall t \in ]t_0-\tau; t_0+\tau[.$ \\
Moreover,
 for all large $N$, \begin{equation}\label{1} ]t_0 -\tau;t_0 + \tau[\cap U_N=]t_0 -\tau;t_0 (N)[.\end{equation}
Since according to Theorem \ref{theoBiane},  $\Psi_N$ is strictly increasing on $U_N$, we have 
\begin{equation}\label{2}\Psi_N(]t_0 -\tau;t_0 (N)[)= ]\Psi_N(t_0 -\tau);\Psi_N(t_0 (N))[.\end{equation} Moreover according to Lemma \ref{croissancePsi},
\begin{equation}\label{3}\Psi_N([t_0 (N);t_0 +\tau [)= [\Psi_N(t_0 (N));\Psi_N(t_0 +\tau)[.\end{equation}
Note that $\Psi_N(t_0 (N))=H_N(t_0(N))=t_0(N) +\int \frac{d\mu_{A_N}(x)}{t_0(N)-x}$
with for $\tau$ small enough and $N$ large enough $t_0(N) \in [t_0-\tau;t_0 +\tau)] \subset \R \setminus \{ \text{supp}(\nu)\cup \Theta\}. $
Lemma \ref{montel} readily yields that $u_0(N)=\Psi_N(t_0 (N))$ converges towards $H(t_0)=\Psi(\Psi^{-1}(u_0))=u_0.$
Now, for $\tau$ small enough, $t_0+\tau \in  \R \setminus \{\overline{U}\cup \Theta\}$ and $t_0 -\tau \in U$, $t_0 -\tau \in  \R \setminus \{ \text{supp}(\nu)\cup \Theta\}$, so that using Lemma \ref{prelim2},
for any $\eta>0$ small enough,  for all large $N$, \begin{equation}\label{4}\Psi_N(t_0+\tau)> u_0 +\eta~~\text{ and }~~\Psi_N(t_0 -\tau)< u_0 -\eta.\end{equation}
It readily follows from  (\ref{1}),  (\ref{2}), (\ref{3}), (\ref{4}), (\ref{densite}) and (\ref{defU}) that for any $\eta>0$ small enough,  for all large $N$, 
$$\forall u  \in 
  [u_0(N); u_0+\eta [ , ~~p_N(u)=0 ~~\text{and} ~~\forall u  \in ]u_0-\eta; u_0(N)[, ~~p_N(u)>0.$$
Now,  for $t$ in a small neighborhood of $t_0$ and $N$ large enough let us define
$$\epsilon_N(t)= \frac{N-r}{N}\int \frac{d\hat{\nu}_N(x)}{(t-x)}-\int \frac{d\nu(x)}{(t-x)}.$$
 We have  \begin{eqnarray*} \Psi_N(t_0 (N))&=&H_N(t_0(N))\\&=&t_0(N) +\int \frac{d\mu_{A_N}(x)}{t_0(N)-x}\\
&=&H(t_0) +f_N(t_0(N)) + \int \frac{d\mu_{A_N}(x)}{t_0(N)-x} -\int \frac{d\nu(x)}{t_0-x}\\&=& H(t_0) +f_N(t_0(N)) + \epsilon_N(t_0(N)) +  \int \frac{d\nu(x)}{t_0(N)-x} -\int \frac{d\nu(x)}{t_0-x}\\&& + O(\frac{1}{N})\\
&=& H(t_0) +f_N(t_0(N)) + \epsilon_N(t_0(N))\\ &&-f_N(t_0(N))\left[ \int \frac{d\nu(x)}{(t_0-x)^2}- f_N(t_0(N)) \int \frac{d\nu(x)}{(t_0(N)-x)(t_0-x)^2} \right]\\ & & + O(\frac{1}{N})\\
&=& H(t_0)  + \epsilon_N(t_0(N))+f_N(t_0(N))^2 \int \frac{d\nu(x)}{(t_0(N)-x)(t_0-x)^2} \\&&+ O(\frac{1}{N})\\&=&H(t_0)  + \epsilon_N(t_0(N)) + \frac{1}{4}({\epsilon_N}^{'}(t_0(N)))^2 (1+o(1))+ O(\frac{1}{N}).
\end{eqnarray*}
The proof of  Proposition \ref{rightedge} is complete
$~~~~~~~\Box$

The proof of Proposition \ref{leftedge} is similar and left to the reader.
\paragraph{Proof of Proposition \ref{approxPearcey}:}
According to Theorem \ref{theoBiane},$$
\Psi^{-1}( ]u_0 -\epsilon; u_0+ \epsilon[) \subset \overline{U} $$
and more precisely, since 
$$p(\Psi (\Psi^{-1}(u)))=\frac{v(\Psi^{-1}(u))}{\pi },$$
we have $\Psi^{-1}( ]u_0 -\epsilon; u_0+ \epsilon[\setminus \{u_0\}) \subset {U} $ and $x_0 =\Psi^{-1}(u_0) \notin U$.
Since we assume that  $\text{supp}(\nu) \subset U$, $x_0 \notin \text{supp}(\nu)$.
 Note that $u_0 =\Psi(x_0) =H(x_0)$ since $v(x_0)=0$.
Moreover, since the homeomorphism $\Psi$ is strictly increasing on $U$, it is easy to see that $\Psi^{-1}$ is  strictly increasing on $]u_0 -\epsilon; u_0+ \epsilon[$ and $\Psi^{-1}( ]u_0 -\epsilon; u_0+ \epsilon[\setminus \{u_0\})=]\Psi^{-1}( u_0 -\epsilon);x_0[\cup ]x_0; \Psi^{-1}( u_0+ \epsilon)[
$. Therefore 
$x_0$ is a point in the complement of $\text{supp}(\nu)$ where two components of the set $U$ merge into one. Therefore, Lemma \ref{theoBiane} implies that $\int \frac{d\nu(s)}{(x_0-s)^2}=1$ and $\int \frac{d\nu(s)}{(x_0-s)^3}=0$.
Since we assume that for any $\theta_i \in \Theta$, $\theta_i \neq x_0$, we have $x_0 \notin \Theta$. Therefore $x_0$ satisfies the assumptions of Lemma \ref{prelim22}. 
Let $\eta$ be such that $0< 2\eta< \epsilon$ and   $[\Psi^{-1}(u_0-2\eta); \Psi^{-1} (u_0+2\eta)] \subset \R \setminus \{\text{supp}(\nu) \cup \Theta\}$.
According to Lemma \ref{prelim22}, for all large $N$,
 there exists $x_0(N) $ in  $
[\Psi^{-1}(u_0-2\eta); \Psi^{-1} (u_0+2\eta)] $ such that $\int \frac{d{\mu_{A_N}}(s)}{(x_0(N)-s)^2}=1$,  $\int \frac{d{\mu_{A_N}}(s)}{(x_0(N)-s)^3}=0$ and  $[\Psi^{-1}(u_0-2\eta); \Psi^{-1} (u_0+2\eta)]  \setminus \{x_0(N)\} \subset U_N$ .
We have also $\lim_{N \rightarrow + \infty}x_0(N)=x_0. $ Note that  since 
$$ \text{for any}~~ x \in \R, ~~p_{N }(\Psi _{N}(x))=\frac{v_{N }(x)}{\pi },$$
we have $$p_{N}(\Psi _{N }(x_0(N)))=0$$ and $$\forall x \in [\Psi^{-1}(u_0-2\eta); \Psi^{-1} (u_0+2\eta)] \setminus \{x_0(N)\}, ~~p_{N }(\Psi _{N }(x))>0.$$
Using Lemma \ref{prelim2}, we can deduce that for all large $N$, $$\Psi_N(\Psi^{-1}(u_0 -2 \eta))< u_0 -\eta ~~\text{and}~~\Psi_N(\Psi^{-1}(u_0 +2 \eta))> u_0 +\eta.$$
Moreover since  $\lim_{N \rightarrow + \infty}x_0(N)=x_0$, we have for all large $N$, $x_0(N) \in ]\Psi^{-1}(u_0-\eta/2); \Psi^{-1} (u_0+\eta/2)[$ so that $u_0(N)= \Psi_N(x_0(N)) \in ]u_0-\eta;u_0+ \eta[$ for all large $N$ by using oncemore Lemma \ref{prelim2}.
The proof is complete $\Box$

%

\paragraph{ Proof of Proposition \ref{outlier}:} According to Lemma \ref{spike},  for $\epsilon>0$ small enough, for all large $N$,  $U_N \bigcap ]\theta_i-\epsilon; \theta_i + \epsilon [= ]t_1^i(N),t_2^i(N)[$ 
where $t_1^i(N)$ and $t_2^i(N)$ satisfy
$$t_1^i(N)= \theta_i  - \sqrt{\frac{k_i}{N} \phi_N(t_1^i(N))}$$
$$t_2^i(N)= \theta_i  + \sqrt{\frac{k_i}{N} \phi_N(t_2^i(N))}$$
with $\phi_N(t)= \frac{1}{1 - \frac{N-r}{N} \int \frac{1}{(t-x)^2}d\hat{\nu}_N (x)- \frac{1}{N} 
 \sum_{j\neq i} \frac{k_j}{( t-\theta_j)^2}}$
and $1 \leq \phi_N(t) \leq K(\epsilon)$ for any $t \in ]\theta_i-\epsilon; \theta_i + \epsilon [$ .
For $N$ large enough $t_1^i(N) > \theta_i -\epsilon/2$ and $t_2^i(N) < \theta_i +\epsilon/2$ and $[\theta_i -\epsilon/2;\theta_i+\epsilon/2]\cap \overline{U_N}=[t_1^i(N),t_2^i(N)]$. Therefore, according to Theorem \ref{theoBiane},
$[\Psi_N(t_1^i(N)),\Psi_N(t_2^i(N))]$ is a connected component of $\text{supp}(\mu_{sc}\boxplus\mu_{A_N})$ and $[\Psi_N(\theta_i -\epsilon/2); \Psi_N(t_1^i(N))[\cup]\Psi_N(t_2^i(N));\Psi_N(\theta_i+\epsilon/2])\subset \R \setminus \text{supp}(\mu_{sc}\boxplus\mu_{A_N})$.
Now, we have
\begin{eqnarray*}
\Psi_N(t_1^i(N)) &= &\theta_i  - \sqrt{\frac{k_i}{N} \phi_N(t_1^i(N))} +
 \frac{N-r}{N}\int \frac{d\hat{\nu}_N(x)}{(\theta_i  - \sqrt{\frac{k_i}{N} \phi_N(t_1^i(N))}-x)} \\
&& + \sum_{j\neq i}
\frac{k_j}{N \left(\theta_i - \theta_j - \sqrt{\frac{k_i}{N} \phi_N(t_1^i(N))}\right)}-
 \frac{k_i}{N  \sqrt{\frac{k_i}{N} \phi_N(t_1^i(N))}}\\
&=& {\theta_i} -\sqrt{\frac{k_i}{N} \phi_N(t_1^i(N))}-  \frac{\sqrt{k_i}}{\sqrt{N}} \frac{1}{\sqrt{\phi_N(t_1^i(N))}}\\
&&+
 \frac{N-r}{N}\int \frac{1}{(\theta_i  - \sqrt{\frac{k_i}{N} \phi_N(t_1^i(N))}-x)} d\hat{\nu}_N (x)
+ O(\frac{1}{N})
\\ 
&=& {\theta_i} - \sqrt{\frac{k_i}{N} \phi_N(t_1^i(N))}\left\{1- \int \frac{1}{(\theta_i-x)^2}d\nu(x)\right\}\\&&
+
  \frac{N-r}{N}\int \frac{1}{\theta_i  -x } d\hat{\nu}_N (x)-  \frac{\sqrt{k_i}}{\sqrt{N}} \frac{1}{\sqrt{\phi_N(t_1^i(N))}}
\\&&+\sqrt{\frac{k_i}{N} \phi_N(t_1^i(N))} \left\{
\frac{N-r}{N}\int \frac{d\hat{\nu}_N (x)}{(\theta_i  -x)^2 } 
- \int \frac{d\nu(x)}{(\theta_i-x)^2}\right\} + O(\frac{1}{N})
\\
&=&\rho_N(\theta_i) - \frac{\tau_i}{\sqrt{N}}+ o(\frac{1}{\sqrt{N}})
\end{eqnarray*}
with $\rho_N(\theta_i):=\frac{1}{N} \sum_{y_j \neq \theta_i} \frac{1}{\theta_i -y_j} +\theta_i$
and  $\tau_i=2\sqrt{k_i}\sqrt{1- \int \frac{1}{(\theta_i-x)^2}d\nu(x)}.$
In the same way $$ \Psi_N(t_2^i(N))=\rho_N({\theta_i}) + \frac{\tau_i}{\sqrt{N}}+ o(\frac{1}{\sqrt{N}}).$$
Note that $ \Psi_N(t_1^i(N))$ and  $\Psi_N(t_2^i(N))$ converges towards $\rho_{\theta_i}=\Psi(\theta_i)$. Since for $\epsilon$ small enough, $[\theta_i -\epsilon; \theta_i+\epsilon] \subset  \R \setminus(\overline{U} \cup \Theta)$ (see (\ref{caract})), according to Lemma \ref{prelim2} (3), $\Psi_N(\theta_i -\epsilon/2 )$ and $\Psi_N(\theta_i +\epsilon/2 )$
converge respectively towards $\Psi(\theta_i -\epsilon/2 )$ and $\Psi(\theta_i +\epsilon/2 )$ and, according to Lemma \ref{croissancePsi},
$\Psi(\theta_i -\epsilon/2 )<\Psi(\theta_i -\epsilon/4 )<\Psi(\theta_i  )<\Psi(\theta_i +\epsilon/4 )<\Psi(\theta_i +\epsilon/2 )$. Now, for all large $N$, $\Psi_N(\theta_i -\epsilon/2 )<\Psi(\theta_i -\epsilon/4 ) $ and $\Psi(\theta_i +\epsilon/4 )<\Psi_N(\theta_i +\epsilon/2 )$.
Then, for any $0<\eta< \min \{\Psi(\theta_i +\epsilon/4 )-\Psi(\theta_i ); \Psi(\theta_i  )-\Psi(\theta_i -\epsilon/4)\}$, for all large $N$, we have $\Psi_N(t_1^i(N))> \Psi(\theta_i) -\eta$ and  $\Psi_N(t_2^i(N))< \Psi(\theta_i)+\eta$  whereas $\Psi_N(\theta_i -\epsilon/2 )< \Psi(\theta_i) -\eta$
and $\Psi_N(\theta_i +\epsilon/2 ) >\Psi(\theta_i) +\eta$. Thus $[\Psi_N(t_1^i(N)),\Psi_N(t_2^i(N))]$ is the unique connected component of $\text{supp}(\mu_{sc}\boxplus\mu_{A_N})$ inside $ ] \Psi(\theta_i) -\eta;  \Psi(\theta_i) +\eta[$.
The proof of Proposition \ref{outlier} is complete.
$\Box$
\section{\label{Sec: asanal}Proofs of Theorem \ref{theo: A}, Theorem \ref{theo: S} and Theorem \ref{theo: P} }

\subsection{Correlation functions of the Deformed GUE}
It is known from Johansson \cite{Johansson} (see also \cite{BrezinHikami}) that the joint eigenvalue density induced by the deformed GUE $M_N$ can be explicitely computed. Furthermore it induces a so-called "determinantal random point field". In other words, if one considers a symmetric function $f: \R^m \to \R$, one has that 
\begin{eqnarray*}&&\E \sum_{1\leq i_1<i_2< \cdots <i_m \leq N}f(\lambda_{i_1}, \ldots, \lambda_{i_m})\cr
&& =\int f(x_1, \ldots, x_m) \frac{1}{m!}\det \left( K_N(x_i, x_j)\right)_{i,j=1}^m \prod_{i=1}^m dx_i,\end{eqnarray*}
where $K_N$ is the so-called correlation kernel of the deformed GUE, which has been explicited by \cite{Johansson}. We here state his result. 

\bp \cite{Johansson} \label{Prop: corrkern}
The correlation of the deformed GUE $M_N$ is given by the double complex integral:
\be 
K_N(u,v)=\frac{N}{(2i \pi)^2}\int_{\Gamma}\int_{\gamma}e^{N\frac{(w-v)^2}{2}-N\frac{(z-u)^2}{2}}\frac{1}{w-z}\prod_{i=1}^N \frac{w-y_i}{z-y_i}dw dz,\ee
where $\Gamma$ encircles the poles $y_1, \ldots, y_N$ and $\gamma$ is a line parallel to the $y-$axis not crossing $\Gamma.$
\ep 

At this point, it is worth mentioning that correlation functions and thus local eigenvalue statistics are invariant through conjugation of the correlation kernel. Indeed, one has that 
$$\det \left (K_N(u_i, u_j)\right )_{i,j=1}^m=\det \left( K_N(u_i, u_j) \frac{h(u_i)}{h(u_j)} \right) _{i,j=1}^m,$$
for any non vanishing function $h$.
This fact will be used many times in this article.

Before starting the asymptotic analysis, we list some important facts and notations that are needed hereafter.

Let $u_0$ be given. Assume that both $u$ and $v$ satisfy $|u-u_0|\leq N^{-{\delta}}$ for some $\delta >0$.
Let us set 
$$F_{u_0}(z):=\frac{(z-u_0)^2}{2}+\int_{\R} \ln (z-y) d\nu(y).$$
Note that $F_{u_0}$ is the first order approximation (as $N \to \infty$) of the true exponential term arising in both $z$ and $w$ integrals in the correlation kernel $K_N$. Indeed the true exponential term arising in both integrals is given by 
$$F_{u_0,N}(z):=\frac{(z-u_0)^2}{2}+\frac{1}{N}\sum_{i=1}^N \ln (z-y_i) .$$
We neglect for a while the fake singularity introduced by the logarithm (as $e^{F_{u_0,N}}$ is holomorphic).
By definition, critical points satisfy 
$$F_{u_0,N}'(z)=z-u_0+\frac{1}{N}\sum_{i=1}^N \frac{1}{z-y_i}=0 $$
and one can note that $F_{u_0,N}''=1-\frac{1}{N}\sum_{i=1}^N \frac{1}{(z-y_i)^2}$ does not depend on $u_0$.
It is also convenient for the following to define the curve of critical points of both $F_u$ and $F_{u,N}$.
Let us define 
$$\mathcal{C}=\{ x\pm i v(x), x\in \R\}.$$
One can check that a critical point of $F_u$ with non null imaginary part lies on 
$$\{x \pm i v(x), x \in U\}=\mathcal{C}\cap \{z \in \C, \Im z\not=0\}=\{ z \in \C, \Im z\not=0, \int \dfrac{1}{|z-y|^2} d\nu(y)=1\}.$$
For any $u \in \Psi(U)$, we denote by $z_c^{\pm}(u)$ these two critical points:
$$z_c^{\pm}(u) = \Psi^{-1}(u) \pm i v( \Psi^{-1}(u)).$$ Formula (\ref{densite}) due to Biane shows that 
$|\Im z_c(u)|=\pi p(u).$\\
If instead $F_u$ has no non real critical point, then $u\in \Psi (U^c)$. As a consequence there exists a unique $z_c(u)\in \mathcal{C}\cap \R={U^c}$ such that $F_u'(z_c(u))=0.$ The real numbers $u$ and $z_c(u)$ are then related by the equation  
$$u:= z_c(u)+\int \frac{1}{z_c(u)-y}d\nu(y) \text{ i.e. }z_c(u)=\Psi^{-1}(u).$$
This follows from the fact that
 $\Psi: \R \rightarrow \R $ is one to one. 
In all cases $z_c^{\pm}(u)$,  $z_c(u)$ and $u$ are related by : $$H(z_c^{(\pm)}(u))=u.$$
Similarly we define $$\mathcal{C}_N=\{ x\pm i v_N(x), x\in \R\}.$$
A critical point of $F_{u,N}$ with non zero imaginary part lies on
$$\{x \pm i v_N(x), x \in U_N\}=\mathcal{C}_N\cap \{z \in \C, \Im z\not=0\}=\{ z \in \C, \Im z\not=0, \frac{1}{N}\sum_{j=1 }^N\frac{1}{|z-y_j|^2}=1\}.$$
For any $u \in \Psi_N(U_N)$, denote by $z_{c,N}^{\pm}(u)$ these two critical points of $F_{u,N}$: $$z_{c,N}^{\pm}(u) = \Psi_N^{-1}(u) \pm i v_N( \Psi_N^{-1}(u)).$$ We note that $F_{u,N}$ necessarily admits $N-1$ other critical points, which are real interlaced with the $y_i$'s. We disregard these critical points. Then one has that
$$H_N(z_{c,N}^{\pm}(u))=u.$$
If instead $F_{u,N}$ has no non real critical points, $u \in \Psi_N(U_N^c)$ and there exists a unique $z_{c,N}(u)\in \R\cap \mathcal{C}_N {=U_N^c}$ such that $F_{u,N}'(z_{c,N}(u))=0.$ Again one has that 
$$u=H_N(z_{c,N}(u))= z_{c,N}(u)+\frac{1}{N}\sum_{i=1}^N \frac{1}{z_{c,N}(u)-y_i}.$$
We emphasize that according to (\ref{caract})
$$\forall z \in  \mathring{ (\mathcal{C}_N \cap \R )}=\overline{U_N}^c,\: \frac{1}{N}\sum_{i=1}^N \frac{1}{(z-y_i)^2}<1,$$
and that, according to  Theorem \ref{theoBiane} and Lemma \ref{croissancePsi},  $u \mapsto \Re z_{c,N}^{(\pm)}(u)=\Psi_N^{-1}(u)$ is a strictly increasing function.

Actually in all the cases we study, it turns out that the critical points, that we here denote by $\mathbf{z_c},$ lie on the real axis. We may therefore need to modify $F_{u,N}$ so that there is no singularity in the logarithm.
It may happen in particular that $\exists\: 1\leq i\leq N$, $y_i<\mathbf{z_c}<y_{i+1}$. However by the assumptions we have made, in all cases there exists $\epsilon >0$ such that $[\mathbf{z_c}-\epsilon,\mathbf{z_c}+\epsilon]  $ contains no eigenvalue $y_j, j=1, \ldots, N$. 
In that case we set 
\be \label{singu}F_{u,N}=\frac{(z-u_0)^2}{2}+\frac{1}{N}\sum_{i: y_i < \mathbf{z_c}+\epsilon} \ln (z-y_i) +\frac{1}{N}\sum_{i: y_i > \mathbf{z_c}+\epsilon} \ln (y_i-z).\ee
The contour $\Gamma$ will be split into two parts: $\Gamma_1$ lying to the left of $\mathbf{z_c}+\epsilon$ and 
$\Gamma_2$ to its right (encircling all the eigenvalues $y_i >\mathbf{z_c}+\epsilon$).
The contour $\gamma$ will be chosen so that it lies to the left of $\mathbf{z_c}+\epsilon$.
All these contours cross the real axis at a point where $F_{u,N}$ has no singularity.
Note that with this new definition of $F_{u,N}$, it is still true that 
$$F_{u,N}'(z)=z-u+\frac{1}{N}\sum_{i=1}^N \frac{1}{z-y_i}.$$ Thus all the subsequent derivatives and the curve $\mathcal{C}_N$ are unchanged with this new definition.
The asymptotic exponential term at $\mathbf{z_c}$ is then given by 
$$F_{u_0}(z)=\frac{(z-u_0)^2}{2}+\int_{(-\infty, \mathbf{z_c}+\epsilon)} \ln (z-y) d\nu(y)+\int_{( \mathbf{z_c}+\epsilon, +\infty)} \ln (-z+y) d\nu(y).$$

\subsection{Asymptotics of the correlation kernel at the edges of the support }

\subsubsection{Proof of Theorem \ref{theo: A}\label{subsec: theoA}}

We start from a right extremity point $d$ of a connected component of $\text{supp}
(\nu \boxplus \mu_{sc})$ so that $p(x)=0, \forall x \in [d, d+\epsilon ]$ for some small $\epsilon >0.$ We assume moreover that for any $\theta_j$ such that $\int \frac{d\nu(s)}{(\theta_i -s)^2} =1$, we have $d \neq \theta_j + m_\nu (\theta_j)$.
According to Proposition \ref{rightedge}, such a point $d$ satisfies $d=H(\mathbf{z_0})$ where $\mathbf{z_0}$ is a real solution of 
$$F_{d}''(\mathbf{z_0})=0.$$ Since $\mathbf{z_0}\notin \text{supp}(\nu)\cup \Theta$, $ (H_3)$ implies that for all large $N$, one also has that $\inf_{k=1, \ldots, N}\text{dist}(\mathbf{z_0},y_k)>0$.
By Proposition \ref{rightedge}, there exists a unique extremity point $d_N$ which is the right endpoint of a connected component of $\text{supp} (\mu_N \boxplus \mu_{sc})$ and such that $|d-d_N|\leq \epsilon$ for any $\epsilon.$
Then there exists a point $\mathbf{z_N}$ such that $$H_N(\mathbf{z_N})=d_N.$$
Let $F_{d_N,N }$ be defined as in (\ref{singu}) with $\mathbf{z_c}=\mathbf{z_N}.$ 
By definition, one has that  $\mathbf{z_N}$ is the real degenerate critical point  associated  to $d_N:$
\begin{equation}
F_{d_N,N }' (\mathbf{z_N})=0, \text{ and }
F_{d_N,N }'' (\mathbf{z_N})=0.
\end{equation}
We now turn to the asymptotics of the correlation kernel.
Let $\alpha \in \R$ to be fixed later. Assume that  
\begin{eqnarray}&u_0:=d_N, &u=u_0+\frac{\alpha x}{N^{\frac{2}{3}}};\: v=u_0+\frac{\alpha y}{N^{\frac{2}{3}}}
 \end{eqnarray}
We assume that there exists a real number $M_0>0$ such that $x, y\geq -M_0.$ If $u_0$ is not the top edge of the support $\text{supp}(\mu_{A_N}\boxplus \mu_{\sigma})$, then $x$ and $y$ shall be bounded from above by $\epsilon_0 N^{2/3}$ with $\epsilon_0$ small enough so that $u_0+\frac{\alpha x}{N^{2/3}}$ is smaller than  the left edge of the next connected component of $\text{supp}(\mu_{A_N}\boxplus \mu_{\sigma})$.\\
The associated rescaled correlation kernel is then $$\frac{\alpha }{N^{\frac{2}{3}}}K_N(u,v).$$
We now consider the asymptotics of the correlation kernel and prove that the rescaled kernel $\frac{\alpha}{N^{2/3}}K_N(u,v)$ uniformly converges to the Airy kernel when $-M_0 \leq x,y \leq \epsilon_0 N^{2/3}.$

\paragraph{}Theorem \ref{theo: A} is an easy consequence of the following Proposition.
Set $$\alpha= 2^{1/3}\frac{1}{|F_{u_0, N}^{(3)}(\mathbf{z_N})|^{1/3}}.$$
$\alpha$ is well defined using Lemma \ref{prelim}, (ii).
\bp \label{prop:airy} There exist constants $q, C, c >0$ such that for any $x,y \in [-M_0, \epsilon_0 N^{2/3}],$
$$\Big | \dfrac{\alpha }{N^{\frac{2}{3}}}K_N(u,v)e^{q(y-x)N^{\frac{1}{3}}} - \mathbf{A}(x,y) \Big | \leq \frac{Ce^{-c(x+y)}}{N^{\frac{1}{3}}},$$
where $ \mathbf{A}$ denotes the Airy kernel.
\ep 
 
\paragraph{Proof of proposition \ref{prop:airy}}
By Cauchy's theory and using the fact  proved in Lemma \ref{prelim} that $$F_{d}^{(3)}(\mathbf{z_0}) =2\int \frac{1}{(\mathbf{z_0}-y)^3}d\nu (y)=a_i>0,$$
one deduces that $F_{u_0, N}^{(3)}(\mathbf{z_N})\geq a_i/2$ and that there exist $a>0, M>0$ and a small $\delta$-neighborhood of $\mathbf{z_N}$
such that 
\begin{equation}\label{above} \forall z, |z-\mathbf{z_N}|\leq \delta, \: \Re F_{u_0, N}^{(3)}(z) > a \text{ and } \Big | F_{u_0, N}^{(4)}(z)\Big | \leq M.\end{equation}

We now rewrite the correlation kernel. To this aim, we split $\Gamma$ into two contours lying respectively to the left and to the right of $\mathbf{z_N}.$ This is possible as we assume that $\Delta:=\inf_{k=1, \ldots, N}\text{dist}(\mathbf{z_0},y_k)>0$ and $|z_N-z_0|<\Delta/2$ for $N$ large enough.
Denote by $\Gamma_1$ the part of the contour $\Gamma$ lying to the left of $\mathbf{z_N}$ and set $\Gamma_2:=\Gamma \setminus \Gamma_1.$
In the correlation kernel given by Proposition \ref{Prop: corrkern}, along $\Gamma_1$, we first rewrite the singularity $$1/(w-z)=\alpha N^{\frac{1}{3}}\int_{\R^+}e^{-N^{\frac{1}{3}}\alpha t_o (w-z)}dt_o,$$ which is valid provided the contour $\gamma$ remains to the right of $\Gamma_1.$
This then yields the following expression for the correlation kernel (up to a conjugation factor):
\begin{eqnarray} 
&\dfrac{\alpha }{N^{\frac{2}{3}}}K_N(u,v)&=\frac{\alpha^2 N^{2/3}}{(2i \pi)^2}\int_{\R^+}dt_o \int_{\Gamma_1}dz \int_{\gamma}dw \cr
&&\quad  e^{-N^{\frac{1}{3}}\alpha t_o (w-z)}e^{N(\frac{w^2}{2}-wv-\frac{z^2}{2}+uz)}\prod_{i=1}^N \frac{w-y_i}{z-y_i}\label{KN1}\\
&&+\frac{\alpha N^{\frac{1}{3}}}{(2i \pi)^2}\int_{\Gamma_2}\int_{\gamma}e^{N\frac{(w-v)^2}{2}-N\frac{(z-u)^2}{2}}\frac{1}{w-z}\prod_{i=1}^N \frac{w-y_i}{z-y_i}dw dz.\cr &&\label{KN2}\end{eqnarray}
We denote by $K_N^{(l)}$ (resp. $K_N^{(r)}(u,v)$) the kernel arising in (\ref{KN1}) (resp. (\ref{KN2}) that we consider separately).\\
Note that it is enough to concentrate on $F_{u_0,N}$ for the saddle point analysis of the correlation kernel. Assume given $q\in \R$ that we will fix later.
We rewrite the correlation kernel (and use conjugation thanks to $q$) as:
\begin{eqnarray}
&&K_N^{(l)}(u,v)e^{q(y-x)N^{\frac{1}{3}}}\cr
&&=\frac{1}{(2i \pi)^2}\int_{\R^+}dt_o \int_{\Gamma_1}\int_{\gamma} H(w, y+t_0) G(z, x+t_0) dw dz,\end{eqnarray}
where
\begin{eqnarray}
&&  H(w, y):=\alpha N^{\frac{1}{3}}e^{N F_{u_0, N}(w)-\alpha y(w-q) N^{\frac{1}{3}}}, \cr
&& G(z, x):=\alpha N^{\frac{1}{3}}e^{-N F_{u_0,N}(z)+\alpha x (z-q) N^{\frac{1}{3}}}.
\end{eqnarray}

Let us first consider the leading term in the exponential defining $H$ and $G$ that is $F_{u_0, N}$.
By the choice of $u_0$, the two first derivatives of the exponential term vanish at the real point $\mathbf{z_N}$ so that standard saddle point analysis suggest that the ascent and descent contours shall be given by lines with direction $(2)i \pi/3$ through the critical point $\mathbf{z_N}$. This is true in a compact neighborhood of $\mathbf{z_N}$, as we see below. We ignore for a while the constraint that the contours do not cross each other.\\
We first check that $\Gamma_1$ and $\gamma$ shall follow the directions $2i\pi/3$ or $i \pi/3$. To consider the constraint that they do not cross each other, we later modify these contours in a $N^{-1/3}$ neighborhood of $\mathbf{z_N}$.
Using (\ref{above}), there exists $\delta_0>0$ and $a=a(\delta_0)$ such that for any $|s|\leq \delta_0$
\begin{eqnarray}\label{estimeedecroissance}
&&\Re \left( F_{u_0, N}(\mathbf{z_N}+se^{i\pi/3})-F_{u_0, N}(\mathbf{z_N})\right) 
\cr &&=-\Re \left (s^3 \int_0^1 \int_0^1 \int_0^1 dt dx dv F_{u_0, N}^{(3)}(\mathbf{z_N}+stxv e^{i\pi/3}) \right )\cr
&&= -s^3 \int_0^1 \int_0^1 \int_0^1 dt dx dv  \Re \frac{2}{N}\sum_{j=1}^N \frac{1}{(\mathbf{z_N} + stx v e^{i\pi/3}-y_j)^3}<-as^3.\cr
&&\Re \left( F_{u_0, N}(\mathbf{z_N}+se^{i2\pi/3})-F_{u_0, N}(\mathbf{z_N})\right) \cr
&&=\Re \left (s^3 \int_0^1 \int_0^1 \int_0^1 dt dx dv F_{u_0, N}^{(3)}(\mathbf{z_N}+stxv e^{2i\pi/3}) \right )\cr
&&= s^3 \int_0^1 \int_0^1 \int_0^1 dt dx dv  \Re \frac{2}{N}\sum_{j=1}^N \frac{1}{(\mathbf{z_N} + stx v e^{2i\pi/3}-y_j)^3}>as^3.
\end{eqnarray}
One can then complete the $w$-contour by a line parallel to the imaginary axis. Indeed one can choose $\delta_0$ small enough so that $\mathbf{z_N}+ \delta_0e^{i\pi/3}$ lies in the domain where $1>\frac{1}{N }\sum \frac{1}{|z-y_i|^2}.$ Thus there exists a constant $a'>0$ such that
$$\frac{d\Re F_{u_0, N}(\mathbf{z_N}+ \delta_0e^{i\pi/3}+it)}{dt}<-a't , \: t>0.$$
As a consequence $\Re F_{u_0, N}$ still decreases along the contour $t\mapsto \mathbf{z_N}+ \delta_0e^{i\pi/3}+it, t>0.$ This yields the descent path $\gamma$ for the $w$-integral. \\
For the $z$-integral, we complete the contour as follows.\\ 
If $\mathbf{z_N}+ \delta_0e^{2i\pi/3}$ lies above the curve $\mathcal{C}_N$, we complete the contour by lines parallel to the real axis $x\mapsto \mathbf{z_N}+ \delta_0e^{2i\pi/3}+x, x<0$, up to the moment one crosses the curve $\mathcal{C}_N$. 
Then this part of contour remains on the domain $\{z, \frac{1}{N}\sum_{j=1}^N \dfrac{1}{|z-y_j|^2}\leq 1\}$.  Thus, one can check that there exists a constant $a''>0$ such that
$$\frac{d\Re F_{u_0, N}(\mathbf{z_N}+ \delta_0e^{2i\pi/3}-x)}{dx}>a''x.$$
This (part of) line is then an ascent path for $F_{u,N}$. \\
At the moment (if it exists) where the curve $x\mapsto \mathbf{z_N}+ \delta_0e^{2i\pi/3}-x, x<0$ crosses $\mathcal{C}_N$, one follows $\mathcal{C}_N$ to the left direction up to the moment of time where $\Im z\leq \delta_0 \sqrt 3 /2$ and then again follow a line parallel to the real axis. 
Due to the fact that $u \mapsto \Re z_c(u)$ is an increasing function, this part of the contour is also an ascent path.\\
If instead $\mathbf{z_N}+ \delta_0e^{2i\pi/3}$ lies below the curve $\mathcal{C}_N$, we first follow the contour 
$\mathbf{z_N}+ \delta_0e^{2i\pi/3}+it, $ where $t\geq 0$ up to the moment one crosses $\mathcal{C}_N$. One then follows $\mathcal{C}_N$ to the left direction up to the moment of time where $\Im z\leq \delta_0 \sqrt 3 /2$ and then again follow a line parallel to the real axis. It is an easy computation to check that this contour is also an ascent path.

\paragraph{}Because $d_N$  may not be the right edge of the support, we need to complete the $z-$contour $\Gamma_2$ to the right of $\mathbf{z_N}$ too. In this case, define $$\mathbf{z_N'}= \inf \{x \in R, x>\mathbf{z_N}, v_N(x)>0\}.$$
Note that the contour $\mathcal{C}_N\cap \{z \in \C, \Re (z)\geq \mathbf{z_N'}\}$ is made of contours around $y_i'$s. 
Let $Z$ be the first point encountered on $\mathcal{C_N}$ to the right of $\mathbf{z_N}'$ such that $\Im (Z)$ is a local maximum.  The contour $\Gamma_2$ then follows $\mathcal{C}_N\cap \{z \in \C, \Re (Z)\geq\Re (z)\geq\mathbf{z_N'}\}$. Afterwards $\Gamma_2$ follows the highest of the two curves $\{Z+x, x>0\}$ and $\mathcal{C}_N\cap\{ \Re (z)>\Re (Z)\}$. The contour is completed by symmetry with respect to the real axis.
Because $\mathcal{C}_N$ is the curve of critical points, along $\Gamma_2$ which lies above $\mathcal{C}_N$, one has that 
\begin{eqnarray*}&&\forall z \in \Gamma_2\cap \mathcal{C}_N, \exists u>u_0, z=z_{c, N} (u)\text{ and }\Re F_{u_0,N}(z) >\Re F_{u_0,N}(\mathbf{z_N});\cr
&& \forall x>0,  \frac{\partial}{\partial x}\Re F_{u_0,N}(Z+x) >0, 
\end{eqnarray*}
as long as $Z+x$ lies above $\mathcal{C}_N.$
This finishes the definition of the contours, apart from the constraint that the two contours cannot cross each other.

\begin{figure} \begin{center}
\begin{picture}(0,0)%
\includegraphics{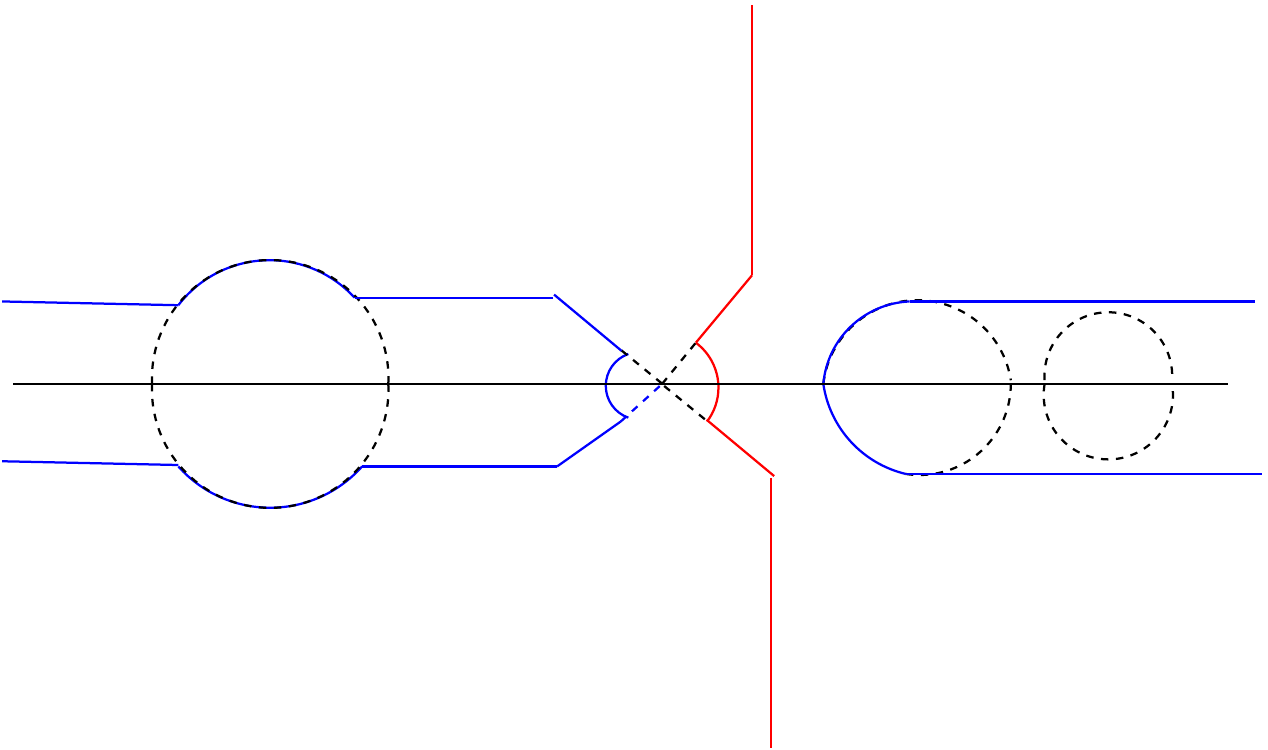}%
\end{picture}%
\setlength{\unitlength}{947sp}%
\begingroup\makeatletter\ifx\SetFigFont\undefined%
\gdef\SetFigFont#1#2#3#4#5{%
  \reset@font\fontsize{#1}{#2pt}%
  \fontfamily{#3}\fontseries{#4}\fontshape{#5}%
  \selectfont}%
\fi\endgroup%
\begin{picture}(25288,14938)(-43,-16680)
\put(18151,-7486){\makebox(0,0)[lb]{\smash{{\SetFigFont{6}{7.2}{\rmdefault}{\mddefault}{\updefault}{\color[rgb]{0,0,1}$\Gamma_2$}%
}}}}
\put(3376,-10486){\makebox(0,0)[lb]{\smash{{\SetFigFont{6}{7.2}{\rmdefault}{\mddefault}{\updefault}{\color[rgb]{0,0,0}$\mathcal{C}_N$}%
}}}}
\put(13051,-10036){\makebox(0,0)[lb]{\smash{{\SetFigFont{6}{7.2}{\rmdefault}{\mddefault}{\updefault}{\color[rgb]{0,0,0}$z_N$}%
}}}}
\put(5251,-6211){\makebox(0,0)[lb]{\smash{{\SetFigFont{6}{7.2}{\rmdefault}{\mddefault}{\updefault}{\color[rgb]{0,0,1}$\Gamma_1$}%
}}}}
\put(15226,-5836){\makebox(0,0)[lb]{\smash{{\SetFigFont{6}{7.2}{\rmdefault}{\mddefault}{\updefault}{\color[rgb]{1,0,0}$\gamma$}%
}}}}
\end{picture}%

\caption{The contour $\Gamma$ and $\gamma$ at an edge.}\label{fig: bord} \end{center}
\end{figure}

\paragraph{}
We now slightly modify the contours in a $N^{-\frac{1}{3}}$ neighborhood of $\mathbf{z_N}$ so that $\gamma$ does not cross $\Gamma_1$. Let $\epsilon>0$ (small) be fixed.  The $w$ and $z$ contours do not go through $\mathbf{z_N}$ but instead follows an arc of circle of ray $\epsilon N^{-\frac{1}{3}}$ centered at $\mathbf{z_N}$ in order to avoid crossing each other (see Figure \ref{fig: bord}).
 We now fix $q=\mathbf{z_N} +\frac{\epsilon}{2} N^{-\frac{1}{3}}$ where $\epsilon$ has been defined as above.
By the estimates on the decay of $F_{u_0,N}$ given in (\ref{estimeedecroissance}), we deduce the following.

Assume first that $|x|, |y| \leq M_0.$ Using (\ref{estimeedecroissance}), we first deduce that there exists $A>0$ such that 
$$\int_{\gamma}H(w, x) dw= \alpha \int_{|w-\mathbf{z_N}|\leq \delta_0}H(w, x)dw (1+O(e^{-AN})).$$
Let us now set $\gamma_0:= \{te^{i\pm \pi/3}, \epsilon \leq t \leq \delta_0 N^{1/3}\} \cup C_{\epsilon}$ where $C_{\epsilon}$ is the arc of circle centered at $0$ joining $\epsilon e^{-i\pi/3}$ and $\epsilon  e^{i\pi/3}$. This contour is oriented from bottom to top. We now make the change of variables $w=\mathbf{z_N}+sN^{-\frac{1}{3}}$
where $s\in \gamma_0.$ We then obtain that 
\begin{eqnarray}
&&\int_{\gamma}H(w, x) dw(1+O(e^{-AN}))\cr
&&=\alpha  \int_{\gamma_0}e^{N F_{u_0,N}(\mathbf{z_N}+sN^{-\frac{1}{3})}-\alpha xs-x\epsilon/2}ds\cr
&&=\alpha \int_{\gamma_0}e^{F_{u_0,N}^{(3)}(\mathbf{z_N})\frac{s^3}{3!} -\alpha xs}e^{NF_{u_0,N}(\mathbf{z_N})-x\epsilon/2 }ds(1+O(N^{-\frac{1}{3}})).
\end{eqnarray}
The last line is obtained by using the fact that 
\begin{eqnarray}&&\Big | e^{N F_{u_0,N}(\mathbf{z_N}+sN^{-\frac{1}{3}})-NF_{u_0,N}(\mathbf{z_N})} - e^{F_{u_0,N}^{(3)}(\mathbf{z_N})\frac{s^3}{3!}}\Big | \cr
&&\leq e^{-as^3} \frac{|s|^4 \sup_{|z-\mathbf{z_N}| \leq \delta_0} |F_{u_0, N}^{(4)}(z)|}{N^{\frac{1}{3}}},
\end{eqnarray} for some constant $a>0$.
More detail can be found in \cite{BBP} Section 3 and we do not develop the computations here.

Similarly we define $\Gamma_0:= \{te^{i2\pm \pi/3}, \epsilon \leq t \leq \delta_0 N^{1/3}\} \cup C'_{\epsilon}$ where $C'_{\epsilon}$ is the arc of circle centered at $0$ joining $ \epsilon e^{-2i\pi/3}$ and $\epsilon e^{2i\pi/3}$. This contour is again oriented from bottom to top.
\begin{eqnarray}
&&\int_{\Gamma_1}G(z, x) dz= \alpha \int_{|z-\mathbf{z_N}|\leq \delta_0}G(z, x)dz (1+O(e^{-AN}))\cr
&&= \alpha \int_{\Gamma_0}e^{-N F_{u_0,N}(\mathbf{z_N}+tN^{-\frac{1}{3})}+\alpha xt+x\epsilon /2}dt(1+O(e^{-AN}))\cr
&&=\alpha \int_{\Gamma_0}e^{-F_{u_0,N}^{(3)}(\mathbf{z_N})\frac{t^3}{3!} +\alpha xt}e^{-NF_{u_0,N}(\mathbf{z_N})+x\epsilon/2}dt(1+O(N^{-\frac{1}{3}})),
\end{eqnarray}
where $t$ describes the contour $\Gamma_0$ formed with the two half lines in the complex plane with angle $e^{\pm 2i \pi/3}$ with respect to the real axis. The contour is also oriented from bottom to top. 
We recall that $\alpha$ has been chosen as
$$\alpha= 2^{1/3}\frac{1}{|F_{u_0, N}^{(3)}(\mathbf{z_N})|^{1/3}}.$$
We then deduce that for $|x|, |y|\leq M_0$, one has that 
\begin{eqnarray*}
&&\Big | \frac{1}{2i \pi}\int_{\gamma}H(w, y)dw -Ai(y)e^{-y\epsilon/2}\Big |\leq \frac{C}{N^{\frac{1}{3}}},\cr
&&\Big | \frac{1}{2i \pi}\int_{\Gamma_1}G(z, y)dz-Ai(y)e^{y\epsilon/2}\Big |\leq \frac{C}{N^{\frac{1}{3}}}.\end{eqnarray*}

We can now conclude to the asymptotic behavior of the rescaled correlation kernel $K_N^{(l)}(u,v)e^{q(y-x)N^{\frac{1}{3}}}$ when $x$ and/or $y$ are allowed to grow unboundedly positive. Indeed for this part of the kernel we do not need to bound $x$ and $y$ from above by $\epsilon_0 N^{2/3}$. 
As by construction the two contours $\Gamma_1$ and $\gamma$ lie respectively to the left (resp. right) strictly of $q$, one can deduce (copying the arguments developed in \cite{BBP} Section 3 ) that 
there exist constants $C, c>0$ such that
\begin{eqnarray}
&&\Big | \frac{1}{2i \pi}\int_{\gamma}H(w, y)dw -Ai(y)e^{-y\epsilon/2}\Big |\leq \frac{C}{N^{\frac{1}{3}}}e^{-c y},\cr
&&\Big | \frac{1}{2i \pi}\int_{\Gamma_1}G(z, y)dz-Ai(y)e^{y\epsilon/2}\Big |\leq \frac{C}{N^{\frac{1}{3}}}e^{-c y}.\label{eq40}\end{eqnarray}
Note that (\ref{eq40}) also holds true (modifying the constants $C,c$ if needed) when $|x|, |y| \leq M_0$.

Last we need to consider the contribution of the contour $\Gamma_2\cup \gamma$. We show that this contribution is negligible provided 
$x$ and $y$ are bounded from above by $\epsilon_0 N^{2/3}$ for some $\epsilon_0$ small enough.
Let us recall that $ \Re F_{u}''(z)\geq 0$ for any $z$ outside $D_N$. Furthermore there exists $\eta>0$ such that  $\text{dist}(\Gamma_2, \gamma)>\eta.$ As a consequence the main contribution from $\Gamma_2$ comes from the closest point to $\mathbf{z_N}$, namely $\mathbf{z_N'}$.   From this we deduce that
\begin{eqnarray}
&&\Big |\alpha \frac{N^{\frac{1}{3}}}{(2i \pi)^2}\int_{\Gamma_2}\int_{\gamma}e^{N(w-v)^2/2-N(z-u)^2/2}\frac{e^{q(y-x)N^{\frac{1}{3}}}}{w-z}\prod_{i=1}^N \frac{w-y_i}{z-y_i}dw dz\Big |\cr
&& \leq C e^{NF_{u_0,N}(\mathbf{z_N})-NF_{u_0,N}(\mathbf{z_N}') +q(y-x)N^{\frac{1}{3}}},\label{eq41}
\end{eqnarray}
for some constant $C>0$.
As $|y|, |x|\leq \epsilon_0 N^{2/3}$, we choose $\epsilon_0>0$ small enough so that there exists a constant $C'>0$ so that
\be \label{eq42}\Re \left ( NF_{u_0,N}(\mathbf{z_N})-NF_{u_0,N}(\mathbf{z_N}') +q(y-x)N^{\frac{1}{3}}\right )<-C'N.\ee
Combining (\ref{eq41}), (\ref{eq42}) and (\ref{eq40}) then yields Proposition \ref{prop:airy}. $\square$

\subsection{Proof of Theorem \ref{theo: S}}

Consider a spike $\theta_{i_1}$ of multiplicity $k_{i_1}$ such that $\int \frac{1}{(\theta_{i_1}-y)^2}d\nu(y)<1$. Then $\theta_{i_1}$ makes $k_{i_1}$ outliers separate from the bulk at $\rho(\theta_{i_1})$ asymptotically, with $\rho(z):=z+\int \frac{1}{z-y}d\nu(y)$. We recall that $\theta_{i_1}$ is such that $\text{dist}(\rho(\theta_{i_1}), \text{supp}(\mu_{sc}\boxplus \nu))>0.$
Thus there exist  (possibly) $\mathbf{z_{N}}$ and $\mathbf{w_{N}}$ such that $H_N(\mathbf{z_{N}})$ and $H_{N}(\mathbf{w_{N}})$ are respectively the right and left endpoints of  the connected component of $\text{supp}(\mu_{sc}\boxplus \mu_{A_N})$ which is respectively on the left hand side and right hand side of $\rho(\theta_{i_1})$
and we have $\mathbf{z_{N}}<\theta_{i_1}<\mathbf{w_{N}}.$ If there is no connected component of 
$\text{supp}(\mu_{sc}\boxplus \mu_{A_N})$ to the right respectively the left of $\rho(\theta_{i_1})$, we then set $\mathbf{w_{N}}=+\infty,$
respectively $z_N =-\infty.$

We first need some definitions to consider the asymptotic correlation functions close to an outlier.
Let $\rho_N$ be defined in (\ref{defrhoN}).
Let $c>0$ be given (to be defined later). 
We set $$u_0:=\rho_N(\theta_{i_1}), \quad u=u_0+\frac{cx}{\sqrt N}, v=u_0+\frac{cy}{\sqrt N}. $$   
Again we assume that $x,y$ are bounded from below by $-M_0$ for some real number $M_0>0$. On the other side, $x$ and $y$ are not allowed to grow unboundedly. Let $\eta_1>0$ be given (small). We assume that $\eta_0 >0$ is small enough so that 
$$\rho_{\theta_{i_1}}+\eta_0< \rho_N(\theta_i)-\eta_1, \forall i \text{ s. t. }\theta_i>\theta_{i_1}, \quad \rho_{\theta_{i_1}}+\eta_0<H_{N}(\mathbf{w_{N}}).$$ We assume that $x,y \leq \eta_0 N^{1/2}.$
We now consider the asymptotics of the rescaled correlation kernel :
$$ \dfrac{c}{\sqrt N}K_N(u,v).$$
Define \be \label{defGUO}G_{u_0,N}(z):= \frac{z^2}{2}-u_0z +\frac{1}{N}\left( \sum_{j: \:y_j<\theta_{i_1}}\ln (z-y_j)\right )+\frac{1}{N}\left( \sum_{j: \:y_j>\theta_{i_1}}\ln (y_j-z)\right ).\ee
We here set $$c:=\sqrt{G_{u_0,N}'' (\theta_{i_1})}>0.$$

Let $K_H$ be the correlation kernel of a $k_{i_1}\times k_{i_1}$ GUE. We recall that $K_H$ is the Christoffel Darboux kernel of some rescaled Hermite polynomials 
 satisfying the orthogonality relationship
$\int_{-\infty}^\infty p_m(x)p_n(x) e^{-\frac12 x^2} dx =
  \delta_{mn}.$
\bp \label{prop: spike}
There exist constants $q,C,$ and $ C'>0$ such that for $x, y \in [-M_0, \eta_0 N^{1/2}]$
$$\Big | \dfrac{c}{\sqrt N}K_N(u,v)e^{qcN^{\frac{1}{2}}(y-x)}-K_H(x,y)\Big | \leq \frac{Ce^{-C'(x+y)}}{N^{1/2}}.$$
\ep 

\paragraph{Proof of Proposition \ref{prop: spike}: }
We again split the correlation kernel into two parts, by dividing the contour $\Gamma$ into two parts. One contour, denoted by $\Gamma_1$ encircles the eigenvalues $y_i$ such that $y_i \leq \theta_{i_1}$. The other contour $\Gamma_2$ then encircles all the eigenvalues $y_j$ such that $y_j>\theta_{i_1}$. This is possible as we assume that spikes are independent of $N$. Note that $\Gamma_1$ can be chosen so that it lies to the left of $\theta_{i_1}+\eta N^{-1/2}$ for some small $\eta>0.$
Accordingly we define $K_N^{(l)}(u,v)$ and $K_N^{(r)}$ to be the corresponding contributions (from contours lying to the left or to the right of $\theta_{i_1}+\eta N^{-1/2}$) to the correlation kernel. 
 
We first rewrite the singularity in the correlation kernel. Then, provided $\Re( w-z)>0$, one has that 
$$\frac{1}{w-z}=\int_{\R^+}dt_0e^{-N^{\frac{1}{2}}ct_0 (w-z)}cN^{\frac{1}{2}}.$$ 
Thus one can write that 
\begin{eqnarray} 
&\dfrac{c}{\sqrt N}K_N^{(l)}(u,v)=&\frac{c^2N}{(2i \pi)^2}\int_{\R^+}dt_0 \: \int_{\Gamma_1}\int_{\gamma}\prod_{i=1}^N \frac{w-y_i}{z-y_i}\cr
&&e^{N(\frac{w^2}{2}-wv)-N(\frac{z^2}{2}-zu)-N^{\frac{1}{2}}ct_0(w-z)}dw dz,\end{eqnarray}
where $\gamma$ is a line parallel to the $y-$axis not crossing $\Gamma_1.$
We keep the other kernel unchanged:
\begin{eqnarray} 
&\dfrac{c}{\sqrt N}K_N^{(r)}(u,v)=&\frac{cN^{1/2}}{(2i \pi)^2} \: \int_{\Gamma_2}\int_{\gamma}\prod_{i=1}^N \frac{w-y_i}{z-y_i}\cr
&&e^{N(\frac{w^2}{2}-wv)-N(\frac{z^2}{2}-zu)}\frac{1}{w-z}dw dz.\end{eqnarray}

Consider the rescaled correlation kernel $\dfrac{c}{\sqrt N}K_N(u,v)e^{qcN^{\frac{1}{2}}(y-x)}$ for some $q$ to be defined. 
We now set, using the definition of $G_{u_0, N}$ given by (\ref{defGUO}):
\begin{eqnarray}&&H(w,y)=c\sqrt N (\sqrt N)^{k_{i_1}}e^{NG_{u_0, N}(w)-N^{\frac{1}{2}}cy(w-q)} (w-\theta_{i_1})^{k_{i_1}},\cr
&&G(z,x)=c\frac{\sqrt N}{(\sqrt N)^{k_{i_1}}} e^{-NG_{u_0, N}(z)+N^{\frac{1}{2}}cx(z-q)}\times (z-\theta_{i_1})^{-k_{i_1}}.
\end{eqnarray}
Then one has that 
\begin{eqnarray} 
&&\dfrac{c}{\sqrt N}K_N^{(l)}(u,v)e^{qcN^{\frac{1}{2}}(y-x)}\cr
&&=\int_{\R^+}dt_0\int_{\gamma}dw\int_{\Gamma_1}dzH(w,y+t_0) G(z,x+t_0) .
\end{eqnarray}

Note that the measure 
$$\tilde{\nu}_N=\frac{1}{N-k_{i_1}}\sum_{j: \:y_j\not= \theta_{i_1}}\frac{1}{z-y_j}$$ still converges to $\nu$.
Let us define $$\tilde{v}_N: \R \mapsto \R, ~~\tilde{v}_N(x) =\inf\{v \geq 0, \int\frac{d\tilde{\nu}_N(s)}{(x-s)^2 +v^2} > \frac{N}{N-k_{i_1}}\},$$
$$\tilde{U}_N=\{ x \in \R, \tilde{v}_N(x)>0\}$$ and  $$\mathcal{C}'_N =\{x\pm i \tilde{v}_N(x), x \in \R\}.$$
In addition $\theta_{i_1}$ is a critical point of $G_{u_0, N}$, which is the leading term in the exponential term defining both $G$ and $H$.
An easy computation shows that 
$G_{u_0,N}''(\theta_{i_1}) >0$. Furthermore one can check that there exist
$\delta>0$ and constants $c(\delta)>0, M(\delta)>0$ such that 
$$\forall \: z , |z-\theta_{i_1}|\leq \delta, \: \Big |G_{u_0, N}''(z)\Big |\geq c(\delta), \text{ and }\Big |G_{u_0, N}^{(3)}(z)\Big |\leq M(\delta).$$

In order to perform the asymptotic analysis of the correlation kernel, we now choose  $$q=\theta_{i_1}+\frac{\epsilon}{2cN^{\frac{1}{2}}}.$$
We start with the kernel $K_N^{(l)}$. We first consider the asymptotics of the function $H.$ We first consider the case where $|x|, |y|\leq M_0$. The other case will be considered hereafter.
Let $\epsilon >0$ be small. 
Define $\gamma=\theta_{i_1}'+it, t\in \R$ oriented from bottom to top where $\theta'_{i_1}=\theta_{i_1}+\frac{\epsilon}{c}  N^{-\frac{1}{2}}.$ One has that  
$$\frac{d}{dt} \Re \left (G_{u_0N}(\theta'_{i_1}+it)\right )= -t \left (1-\frac{1}{N}\sum_{j: \:y_j\not= \theta_{i_1}}\frac{1}{|\theta'_{i_1}-y_j +it|^2}\right ) \leq -Ct,$$
for some constant $C>0$. This follows from the fact that the second derivative of $G_{u_0,N}$ does not vanish in a neighborhood of $\theta_{i_1}$ in particular. 
Note also that the variation of $G_{u_0,N}(\theta'_{i_1})-G_{u_0,N}(\theta_{i_1})$ is of the order of $1/N$.
We now use the same arguments as in Subsection \ref{subsec: theoA}. As we see just below, we can deform $\Gamma_1$ so that $\gamma$ lies strictly to the right of $\Gamma_1$.
Assuming this holds true, one gets that there exists a constant $A>0$ such that
\begin{eqnarray*}&&\int_{\gamma}H(w,y)dw\cr
&&=  c(\sqrt N)^{k_{i_1}+1}\int_{ |w-\theta'_{i_1}|\leq \delta}e^{NG_{u_0,N}(w)-N^{\frac{1}{2}}cy(w-q)} \times (w-\theta_{i_1})^{k_{i_1}}(1+O(e^{-AN})).\end{eqnarray*}
Making the change of variables $w=\theta_{i_1}+i\frac{t}{c\sqrt N}$,  and setting  $\R_{\text{def}}= \R -i \epsilon$ one obtains that 
\begin{eqnarray*}&&\int_{\gamma}H(w,y)dw(1+O(e^{-AN}))\cr
&&=ce^{N G_{u_0,N}(\theta_{i_1})}e^{y \epsilon/2}\int_{\R_{\text{def}}}\frac{i}{c}e^{-\frac{t^2}{2}-yit}
(i\frac{t}{c})^{k_{i_1}}(1+O(N^{-\frac{1}{2}})) \cr
&&=e^{N G_{u_0, N}(\theta_{i_1})}\int_{\R_{\text{def}}}i e^{-\frac{t^2}{2}-y(it-\epsilon/2)}  \left( \frac{i t}{c}\right)^{k_{i_1}} (1+O(N^{-\frac{1}{2}})).
\end{eqnarray*}

We consider now the case where $y$ can be as large as $\epsilon_0 N^{1/2}$. We use the fact that the contour $\gamma$ remains to the right of $q$ strictly.
In particular, one can show that there exist constants $C, C'>0$ such that 
\be \Big | \int_{\gamma}\frac{H(w,y)}{e^{N G_{u_0,N}(\theta_{i_1})}}dw-\int_{\R_{\text{def}}}ie^{-\frac{t^2}{2}-y(it-\epsilon/2)}  \left( \frac{i t}{c}\right)^{k_{i_1}}\Big |  \leq \frac{Ce^{-C' y}}{\sqrt N}. \label{est1}\ee

\begin{figure} 
\begin{center}

\begin{picture}(0,0)%
\includegraphics{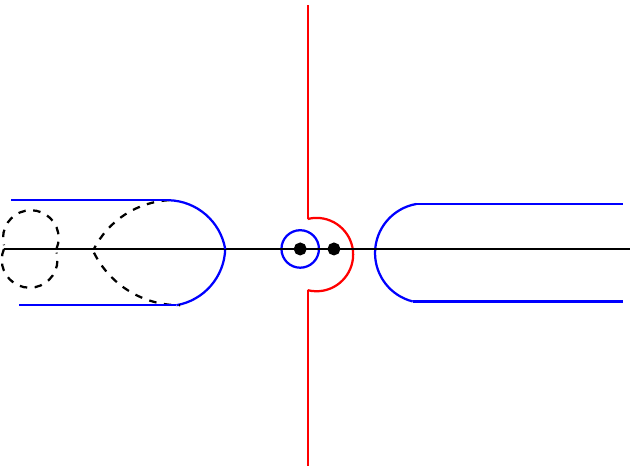} %
\end{picture}%
\setlength{\unitlength}{947sp}%
\begingroup\makeatletter\ifx\SetFigFont\undefined%
\gdef\SetFigFont#1#2#3#4#5{%
  \reset@font\fontsize{#1}{#2pt}%
  \fontfamily{#3}\fontseries{#4}\fontshape{#5}%
  \selectfont}%
\fi\endgroup%
\begin{picture}(12645,9313)(14400,-15180)
\put(20851,-10486){\makebox(0,0)[lb]{\smash{{\SetFigFont{6}{7.2}{\rmdefault}{\mddefault}{\updefault}{\color[rgb]{0,0,0}$q$}%
}}}}
\put(15601,-9211){\makebox(0,0)[lb]{\smash{{\SetFigFont{6}{7.2}{\rmdefault}{\mddefault}{\updefault}{\color[rgb]{0,0,1}$\Gamma_1$}%
}}}}
\put(15676,-11611){\makebox(0,0)[lb]{\smash{{\SetFigFont{6}{7.2}{\rmdefault}{\mddefault}{\updefault}{\color[rgb]{0,0,0}$\mathcal{C}_N$}%
}}}}
\put(20701,-7636){\makebox(0,0)[lb]{\smash{{\SetFigFont{6}{7.2}{\rmdefault}{\mddefault}{\updefault}{\color[rgb]{0,0,1}$\gamma$}%
}}}}
\put(23776,-9511){\makebox(0,0)[lb]{\smash{{\SetFigFont{6}{7.2}{\rmdefault}{\mddefault}{\updefault}{\color[rgb]{0,0,1}$\Gamma_2$}%
}}}}
\put(19651,-11461){\makebox(0,0)[lb]{\smash{{\SetFigFont{6}{7.2}{\rmdefault}{\mddefault}{\updefault}{\color[rgb]{0,0,1}$\theta_{i_1}$}%
}}}}
\end{picture}%
\caption{The contour $\Gamma$ and $\gamma$ at a spike.}\label{fig: spike} 
\end{center}
\end{figure}
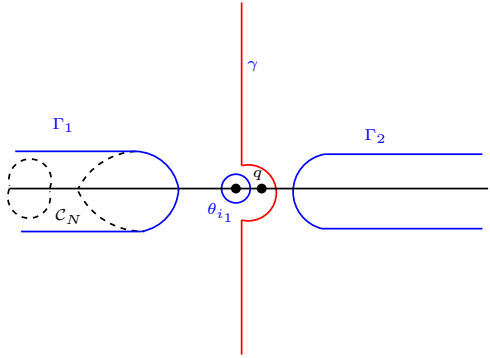

We now turn to the asymptotics of $\int_{\Gamma_1}G(z,y) dz$.
Similarly for the $z$ contour, we use the following contour $\Gamma_1$ (see Figure \ref{fig: spike}).

 First $\Gamma_1$ contains a circle of ray $\frac{\epsilon}{4 cN^{\frac{1}{2}}}$ around $\theta_{i_1}$.
$\Gamma_1$ then has to encircle all the eigenvalues to the left of $\theta_{i_1}$. Note that there exists $\eta>0$ 
$$\sup \{x \in \tilde{U}_N, x <\theta_{i_1}\} =: \mathbf{w_N'}\leq \theta_{i_1}-\eta $$ and $$\inf \{x \in \tilde{U}_N, x >\theta_{i_1}\} =: \mathbf{z_N'} \geq \theta_{i_1}+\eta. $$
Let then $Z'$ be the first point along $\mathcal{C}'_N$ to the left of  $\mathbf{w_N'}$ such that $\Im(Z')$ is a local maximum. 
$\Gamma_1$ then follows $\mathcal{C}'_N$  from  $\mathbf{w_N'}$ to the left direction up to $Z'$.
Then to the left of $Z'$, $\Gamma_1$ follows the highest of the two curves $\mathcal{C}'_N$   and $Z'-x, x>0$. The contour is completed by symmetry with respect to the real axis.
Computing residues, one easily gets that the asymptotics for $G(z, y)$ splits into two parts\\
-the residue at $\theta_{i_1}$ that yields by a straightforward Taylor approximation:
$$e^{-\frac{\epsilon}{2} y}e^{-N G_{u_0,N}(\theta_{i_1})}\text{Res}_{a=0}\left( \big (\frac{c}{a}\big )^{k_{i_1}}e^{-NG_{u_0}(\theta_{i_1}+\frac{a}{c \sqrt N})+NG_{u_0,N}(\theta_{i_1})+ya} \right).$$ 
-The contribution of the rest of the contour $\Gamma_1 \cap \{z, \in \C, \Re z< \theta_{i_1}-\eta\}$ which, by a small extension of the previous subsection, is in the order of 
$$ e^{-NG_{u_0,N}(\mathbf{w_N'})}<< e^{-NG_{u_0,N}(\theta_{i_1})}.$$
This is also exponentially negligible in the large $N$ limit.

\paragraph{}To finish the asymptotic analysis of $G$, we show that the first term is indeed in the order of $e^{-N G_{u_0,N}(\theta_{i_1})}.$
By a straightforward Taylor expansion one obtains that 
\begin{eqnarray}&& e^{-\frac{\epsilon y}{2}}\Big | \text{Res}_{a=0}\left( \big (\frac{c}{a}\big )^{k_{i_1}}e^{N(G_{u_0,N}(\theta_{i_1})-G_{u_0,N}(\theta_{i_1}+\frac{a}{c N^{\frac{1}{2}}}))+ya} \right)-
\text{Res}_{a=0}\left( \big (\frac{c}{a}\big )^{k_{i_1}}e^{ay-\frac{a^2}{2}} \right) \Big | \cr
&&\leq \frac{Ce^{-C' y}}{\sqrt N}, \label{est2}
\end{eqnarray}
for some constants $C, C'>0$. The exponential decay for large $y$ follows again from the fact that the residue is computed on a circle of ray $\epsilon /4c$ lying to the left strictly of $\epsilon /2c$.

We now turn to the asymptotic analysis of $K_N^{(r)}(u,v).$ Let us define the contour $\Gamma_2$ as in the preceding section. Let $Z$ be 
 the first point along $\mathcal{C}'_N$ to the left of $\mathbf{z_N'}$ such that $\Im(Z')$ is a local maximum. 
$\Gamma_2$ first  follows the part $\mathcal{C}'_N$ lying to the right of $\mathbf{z_N'}$ up to the moment where it reaches $Z$. Then $\Gamma_2$ is pursued to the right by following the highest of the two curves $\mathcal{C}'_N$ and $Z+x, x>0$. Again it is completed by symmetry with respect to the real axis. 
It is an easy computation to check that $\Re G_{u_0}(z)$ achieves its minimum on $\Gamma_2$ at $\mathbf{z_N'}$.
The contour $\gamma$ is chosen as before. Note that the function $\frac{1}{w-z}$ remains bounded along $\gamma \cup \Gamma_2.$
We then deduce that 
\begin{eqnarray} 
&\Big |\dfrac{c}{\sqrt N}K_N^{(r)}(u,v)e^{N^{1/2}c(y-x)q}\Big |=&\leq C e^{N \Re G_{u_0,N}(\theta_{i_1})- G_{u_0,N}(\mathbf{z_N'})+N^{1/2}c(y-x)q}\cr &&\leq Ce^{-C'N},\label{lastkern}
\end{eqnarray}
provided $\epsilon_0$ is small enough. 
Thus the kernel $\dfrac{c}{\sqrt N}K_N^{(r)}(u,v)e^{N^{1/2}c(y-x)q}$ converges uniformly to $0$ on $[-M_0, \epsilon_0 N^{1/2}]$.
Combining (\ref{est1}), (\ref{est2}) and (\ref{lastkern}) then yield Proposition \ref{prop: spike} using the expression of the correlation functions of $k_{i_1}\times k_{i_1}$ GUE given in Section 4.3 of \cite{BBP}. $\square$

\subsection{At a point where two connected components merge}
Let now consider a point $u \in \text{supp}(\mu_{\sigma}\boxplus \nu)$ such that the density $p$ of $\mu_{\sigma}\boxplus \nu$ verifies
 $$ p(u) =0, \:p (x) >0\: \forall\: x\in [u-\epsilon/2, u+\epsilon/2]\setminus \{u\} \text{ for some $\epsilon >0$}.$$
This means that the critical point $z_c(u)$ associated to $u=H(z_c(u))$ is unique, real and lies at the "intersection" of two complex curves (see Figure \ref{fig: dessinPear} below).
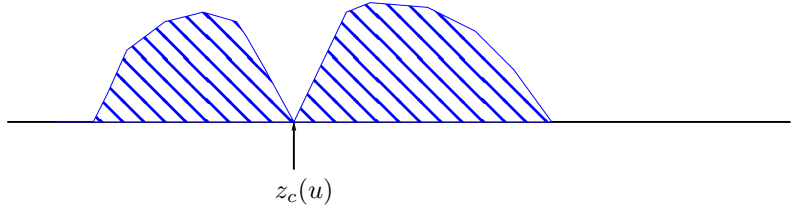
\begin{figure}[h!]\begin{center}
\begin{picture}(0,0)%
\includegraphics{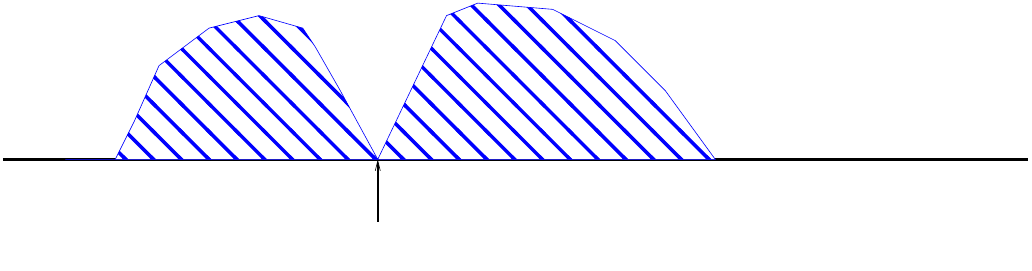}%
\end{picture}%
\setlength{\unitlength}{1579sp}%
\begingroup\makeatletter\ifx\SetFigFont\undefined%
\gdef\SetFigFont#1#2#3#4#5{%
  \reset@font\fontsize{#1}{#2pt}%
  \fontfamily{#3}\fontseries{#4}\fontshape{#5}%
  \selectfont}%
\fi\endgroup%
\begin{picture}(12366,3223)(493,-6497)
\put(4726,-6361){\makebox(0,0)[lb]{\smash{{\SetFigFont{10}{12.0}{\rmdefault}{\mddefault}{\updefault}{\color[rgb]{0,0,0}$z_c(u)$}%
}}}}
\end{picture}%
\caption{ A point in the bulk with vanishing density.}\label{fig: dessinPear}
\end{center}
\end{figure}
Because $z_c(u)\notin \text{ supp}(\nu),$ we deduce from Lemma \ref{lemmeBiane} that 
$$F''(z_c(u))=0 \text{ and that  }F^{(3)}(z_c(u))=0. $$
The first order derivative which does not vanish at $z_c(u)$ is then the fourth one: $F^{(4)}(z_c(u))<0.$
For the asymptotic exponential term $F$, $z_c(u)$ is a doubly degenerate critical point. 
Thanks to Proposition \ref{approxPearcey}, one can transmit this double degeneracy to the true exponential term $F_{u,N}$.
There exists a unique point $z_{c, N}$  in a $\eta$-neighborhood of $z_c$ (for any $\eta>0$) such that 
$$F_{u,N}''(z_{c,N})=F_{u,N}^{(3)}(z_{c,N})=0.$$ At such a point, one obviously has that 
$$F_{u,N}^{(4)}(z_{c,N})<0.$$
Here $F_{u,N}$ is defined by (\ref{singu}) with $\mathbf{z_c}=z_{c,N}.$
Set $u_0=H_N(z_{c, N}).$ We here show that the asymptotic correlation functions in the vicinity of $u_0$ are determined by the so-called Pearcey kernel defined by (\ref{defKP}).
\bp \label{prop: pearcey} Set $\kappa= |F^{(4)}(z_{c,N})|^{1/4} $. Uniformly for $x, y$ in a fixed compact interval, one has that 
$$\lim_{N \to \infty }\frac{\kappa}{N^{\frac{3}{4}}}K_N(u_0+\frac{\kappa x}{N^{\frac{3}{4}}}, u_0+\frac{\kappa y}{N^{\frac{3}{4}}})=K_P(x,y).$$
\ep 

\begin{figure}[h!]\begin{center}
\begin{picture}(0,0)%
\includegraphics{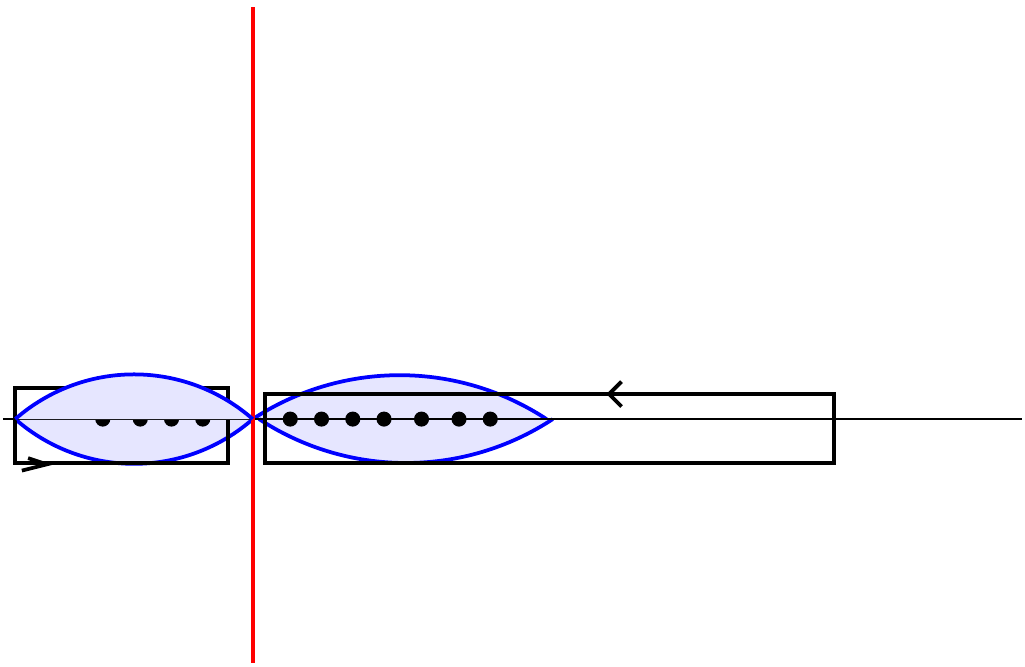}%
\end{picture}%
\setlength{\unitlength}{1579sp}%
\begingroup\makeatletter\ifx\SetFigFont\undefined%
\gdef\SetFigFont#1#2#3#4#5{%
  \reset@font\fontsize{#1}{#2pt}%
  \fontfamily{#3}\fontseries{#4}\fontshape{#5}%
  \selectfont}%
\fi\endgroup%
\begin{picture}(12291,7963)(268,-7905)
\put(2701,-2086){\makebox(0,0)[lb]{\smash{{\SetFigFont{10}{12.0}{\rmdefault}{\mddefault}{\updefault}{\color[rgb]{1,0,0}$\gamma$}%
}}}}
\put(8101,-4261){\makebox(0,0)[lb]{\smash{{\SetFigFont{10}{12.0}{\rmdefault}{\mddefault}{\updefault}{\color[rgb]{0,0,0}$\Gamma$}%
}}}}
\put(4726,-5236){\makebox(0,0)[lb]{\smash{{\SetFigFont{10}{12.0}{\rmdefault}{\mddefault}{\updefault}{\color[rgb]{0,0,0}$y_i$}%
}}}}
\end{picture}%

\caption{Initial contours $\gamma$ and $\Gamma$, which do not cross.}\label{fig: pearcini} 
\end{center}
\end{figure}
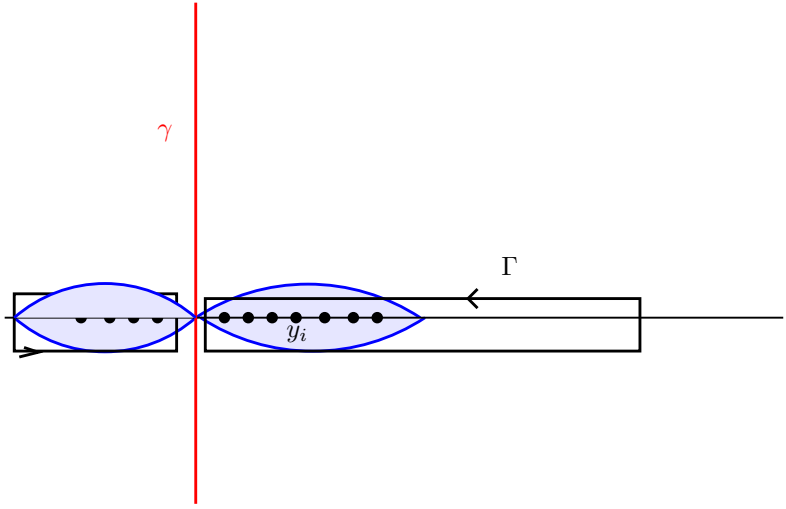

\paragraph{Proof of Proposition \ref{prop: pearcey}}
We start from the expression for the correlation kernel given in Proposition \ref{Prop: corrkern}, where the contours are as shown on Figure \ref{fig: pearcini}.


One has that $F_N^{(4)}(z_{c,N})<0$ and it is not difficult to see that, given $\delta >0$ small, there exists a constant $M$ such that $|F_N^{(5)}(z)|\leq M$ for all complex numbers $z,$ such that  $|z-z_{c,N}|\leq \delta.$
From this we deduce that for any  real $t$ such that $|t|\leq \delta $
$$\Big |F_{u,N}(z_{c,N}+te^{i\frac{\pi}{4}})-F_{u,N}(z_{c,N})+F_{u,N}^{(4)}(z_{c,N})\frac{t^4}{4!} \Big |\leq \frac{M|t|^5}{5!}.$$
Assume that $|t| \leq \delta$,
then one has that $$\Re \left ( F_{u,N}(z_{c,N}+te^{i\frac{\pi}{4}})-F_{u,N}(z_{c,N})\right) \geq |F_N^{(4)}(z_{c,N})|t^4/8!,$$
provided $\delta$ is small enough. 
This ensures that the $(z$-)contour made of two lines with direction $\pm \pi /4$ with the real axis is an ascent contour for $F_{u,N}$, at least in a $\delta$ neighborhood of $z_{c, N}$. To complete the $z$-contour, we need to encircle all the remaining eigenvalues. We pursue the contour as before. 
If $z_{c,N}+\delta e^{i\frac{\pi}{4}}$ (resp. $z_{c,N}+\delta e^{3i\frac{\pi}{4}}$) lies above $\mathcal{C}_N$, the contour goes
 parallely to the real axis to the right (resp. left) up to the moment of time one crosses the curve $\mathcal{C}_N$. 
 Then it follows $\mathcal{C}_N$ to the right (resp. left) direction up to the moment where it crosses the line $\Im z=\delta \sqrt{2}/2$ and so on. 
 If instead $z_{c,N}+\delta e^{i\frac{\pi}{4}}$ (resp. $z_{c,N}+\delta e^{3i\frac{\pi}{4}}$) lies below $\mathcal{C}_N$, then one first joins $\mathcal{C}_N$ along $z_{c,N}+\delta e^{i\frac{\pi}{4}}+it , t\geq 0$ (resp. $z_{c,N}+\delta e^{3i\frac{\pi}{4}}+it, t\geq 0$) and then follows $\mathcal{C}_N$ to the right (resp. left) direction (not going below the line $\Im z=\delta \sqrt{2}/2).$ 
The contour is then completed by symmetry with respect to the real axis. 

For the $w$ contour it is an easy computation that the curve $z_{c,N}+it, t\in \R$ satisfies the descent assumption.
Last, so that the $w$ and $z$ contours do not cross each other, we deform the $z$ contour in a small neighbordhood of $z_{c,N}$ 
to the new contour  $\Gamma_0$ as on Figure \ref{fig: contP}.

We can now conclude to the asymptotic behavior of the kernel. We make the change of variables $w=z_{c,N}+sN^{-1/4}$, $z=z_{c,N}+tN^{-1/4}$, neglecting the part of the contour where $|w-z_{c,N}|\geq \delta$ or $|z-z_{c,N}|\geq \delta$.\\
One has that (up to a conjugation factor)
\begin{eqnarray}
&&\frac{1}{N^{\frac{3}{4}}}K_N(u+\frac{x}{N^{\frac{3}{4}}}, u+\frac{y}{N^{\frac{3}{4}}})\cr &&=\frac{1}{(2 i \pi )^2}\int_{\mathcal{D}}dt \int_{i \R}ds
e^{F^{(4)}(z_{c,N}) \frac{s^4-t^4}{4!} -sy+tx}\frac{1}{s-t}(1+O(N^{-\frac{1}{4}})),
\end{eqnarray}
where we first neglected the parts of the contour lying at a distance $\delta >0$ of $z_{c,N}$ and then performed a Taylor expansion, using the boundedness of the fifth derivative $F_{u,N}^{(5)}$ in a compact neighborhood of $z_{c,N}$. The last estimate holds uniformly for $x,y$ in a fixed compact real interval. Then making the change of variables $s=|F^{(4)}(z_{c,N})|^{1/4} s'$ yields the desired result. $\square$

\end{document}